\documentclass[12pt,reqno]{amsart}
\usepackage{amsthm,amsfonts,amssymb,euscript}

\newcommand{\R}{\ensuremath{\mathbb{R}}}
\newcommand{\Z}{\ensuremath{\mathbb{Z}}}

\newcommand{\e}{\ensuremath{\mathbf{e}}}

\newcommand{\FF}{\ensuremath{\mathcal{F}}}

\numberwithin{equation}{section}
\begin{document}

\title[Global Schr\"{o}dinger maps]
{Global existence and uniqueness of Schr\"{o}dinger maps in dimensions $d\geq 4$}
\author{I. Bejenaru}
\address{University of California -- Los Angeles}
\email{bejenaru@math.ucla.edu}
\author{A. D. Ionescu}
\address{University of Wisconsin -- Madison}
\email{ionescu@math.wisc.edu}
\author{C. E. Kenig}
\address{University of Chicago}
\email{cek@math.uchicago.edu}
\thanks{The second author was
supported in part by an NSF grant and a Packard fellowship.
The third author was supported in part by an NSF grant.}

\begin{abstract}
In dimensions $d\geq 4$, we prove that the Schr\"{o}dinger map initial-value problem
\begin{equation*}
\begin{cases}
&\partial_ts=s\times\Delta s\,\text{ on }\,\mathbb{R}^d\times\mathbb{R};\\
&s(0)=s_0
\end{cases}
\end{equation*}
admits a unique solution $s:\mathbb{R}^d\times\mathbb{R}\to\mathbb{S}^2\hookrightarrow\mathbb{R}^3$, $s\in C(\mathbb{R}:H^{\infty}_Q)$, provided that $s_0\in H^{\infty}_Q$ and $\|s_0-Q\|_{\dot{H}^{d/2}}\ll 1$, where $Q\in\mathbb{S}^2$.

\end{abstract}
\maketitle
\tableofcontents

\section{Introduction}\label{intro}

In this paper we consider the Schr\"{o}dinger map initial-value problem
\begin{equation}\label{Sch1}
\begin{cases}
&\partial_ts=s\times\Delta s\,\text{ on }\,\mathbb{R}^d\times\mathbb{R};\\
&s(0)=s_0,
\end{cases}
\end{equation}
where $d\geq 4$ and
$s:\mathbb{R}^d\times\mathbb{R}\to\mathbb{S}^2\hookrightarrow\mathbb{R}^3$
is a continuous function. The Schr\"{o}dinger map equation has a
rich geometric structure and arises naturally in a number of
different ways; we refer the reader to \cite{NaStUh} or
\cite{KePoStTo} for details. 

For $\sigma\geq 0$ and $n\in\{1,2,\ldots\}$ let $H^{\sigma}=H^\sigma(\mathbb{R}^d;\mathbb{C}^n)$ denote the Banach spaces of $\mathbb{C}^n$-valued Sobolev functions on $\mathbb{R}^d$, i.e.
\begin{equation*}
H^{\sigma}=\{f:\mathbb{R}^d\to\mathbb{C}^n:\|f\|_{H^\sigma}=\big[\sum_{l=1}^n\|\mathcal{F}_{(d)}(f_l)\cdot ( |\xi|^2+1)^{\sigma/2}\|_{L^2}^2\big]^{1/2}<\infty\},
\end{equation*}
where $\FF_{(d)}$ denotes the Fourier
transform on $L^2(\R^d)$. For $\sigma\geq 0$, $n\in\{1,2,\ldots\}$, and $f\in H^\sigma(\mathbb{R}^d;\mathbb{C}^n)$, we define
\begin{equation*}
\begin{split}
\|f\|_{\dot{H}^{\sigma}}=\big[\sum_{l=1}^n\|\FF_{(d)}(f_l)(\xi)\cdot |\xi|^\sigma\|^2_{L^2}\big]^{1/2}.
\end{split}
\end{equation*}
For $\sigma\geq 0$ and $Q=(Q_1,Q_2,Q_3)\in\mathbb{S}^2$ we define the complete metric space
\begin{equation}\label{Sch2}
H^\sigma_Q=H^{\sigma}_Q(\mathbb{R}^d;\mathbb{S}^2\hookrightarrow\mathbb{R}^3)=\{f:\mathbb{R}^d\to\mathbb{R}^3:|f(x)|\equiv 1\text{ and }f-Q\in H^\sigma\},
\end{equation}
with the induced distance
\begin{equation}\label{Ban2}
d^\sigma_Q(f,g)=\|f-g\|_{H^\sigma}.
\end{equation}
For simplicity of notation, we let $\|f\|_{H^\sigma_Q}=d^\sigma_Q(f,Q)$ for $f\in H^\sigma_Q$. Let $\Z_+=\{0,1,\ldots\}$. For $n\in\{1,2,\ldots\}$ and $Q\in\mathbb{S}^2$ we define the complete metric spaces 
\begin{equation*}
H^\infty=H^\infty(\R^d;\mathbb{C}^n)=\bigcap_{\sigma\in\Z_+}H^\sigma\,\,\,\text{ and }\,\,\,H^\infty_Q=\bigcap_{\sigma\in\Z_+}H^\sigma_Q,
\end{equation*}
with the induced distances. Our main theorem concerns global existence and uniqueness of solutions of the initial-value problem \eqref{Sch1} for data $s_0\in H^{\infty}_Q$, with $\|s_0-Q\|_{\dot{H}^{d/2}}\ll1$.

\newtheorem{Main1}{Theorem}[section]
\begin{Main1}\label{Main1} Assume $d\geq 4$ and $Q\in\mathbb{S}^2$. Then there is $\varepsilon_0=\varepsilon_0(d)>0$ such that for any $s_0\in H^{\infty}_Q$ with $\|s_0-Q\|_{\dot{H}^{d/2}}\leq \varepsilon_0$
there is a unique solution
\begin{equation}\label{amy1}
s=S_Q(s_0)\in
C(\R:H^{\infty}_Q)
\end{equation}
of the initial-value problem \eqref{Sch1}. Moreover
\begin{equation}\label{amy2}
\sup_{t\in\R}\|s(t)-Q\|_{\dot{H}^{d/2}}\leq C\|s_0-Q\|_{\dot{H}^{d/2}},
\end{equation}
and 
\begin{equation}\label{amy3}
\sup_{t\in[-T,T]}\|s(t)\|_{H^\sigma_Q}\leq C(\sigma,T,\|s_0\|_{H^\sigma_Q})
\end{equation}
for any $T\in[0,\infty)$ and $\sigma\in\Z_+$.
\end{Main1}

{\bf{Remark:}} We prove in fact a slightly stronger statement: there is $\sigma_0\in[d/2,\infty)\cap\Z$ sufficiently large such that for any $s_0\in H^{\sigma_0}_Q$ with $\|s_0-Q\|_{\dot{H}^{d/2}}\leq \varepsilon_0$
there is a unique solution
\begin{equation*}
s=S_Q(s_0)\in C(\R:H^{\sigma_0-1}_Q)\cap L^\infty(\R:H^{\sigma_0}_Q)
\end{equation*}
of the initial-value problem \eqref{Sch1}. Moreover, the bounds \eqref{amy2} and \eqref{amy3} (assuming $s_0\in H^\sigma_Q$, $\sigma\in\Z_+$) still hold.

The main point of Theorem \ref{Main1} is the global (in time) existence of solutions. Its direct analogue in the setting of wave maps is the work of Tao \cite{Ta1} (see also \cite{
KlMa}, \cite{KlSe}, \cite{Tat1}, \cite{Tat2}, \cite{Ta2}, \cite{KlRo}, \cite{ShSt}, \cite{NaStUh3}, and \cite{Tat3} for other local and global existence (or well-posedness) theorems for wave maps). However, our proof of Theorem \ref{Main1} is closer to that of \cite{ShSt} and \cite{NaStUh3}.

The initial-value problem \eqref{Sch1} has been studied extensively (also in the
 case in which the sphere $\mathbb{S}^2$ is replaced by more general targets). It is known that sufficiently smooth solutions exist locally in time, even for large data (see, for example, \cite{SuSuBa}, \cite{ChShUh}, \cite{DiWa2}, \cite{Ga}, \cite{KePoStTo} and the references therein). Such theorems for (local in time) smooth solutions are proved using delicate geometric variants of the energy method. For low-regularity data, the initial-value problem \eqref{Sch1} has been studied indirectly using the ``modified Schr\"{o}dinger map equations'' (see, for example, \cite{ChShUh}, \cite{NaStUh}, \cite{NaStUh2}, \cite{KeNa}, \cite{Ka}, and \cite{KaKo}) and certain enhanced energy methods. 

In \cite{IoKe2}, Ionescu--Kenig realized that the initial-value problem \eqref{Sch1} can be analyzed perturbatively using the stereographic model, in the case of ``small data'' (i.e. data that takes values in a small neighborhood of a point on the sphere), and proved local well-posedness for small data in $H^\sigma_Q$, $\sigma>(d+1)/2$, $d\geq 2$. The resolution spaces constructed in \cite{IoKe2} (see also \cite{IoKe} for the $1$-dimensional version of these spaces) are based on directional $L^{p,q}_\e$ physical spaces, which are related to local smoothing; in particular, the nonlinear analysis is based on local smoothing and the simple inclusion
\begin{equation*}
L^{\infty,2}_\e\cdot L^{2,\infty}_\e\cdot L^{2,\infty}_\e\subseteq L^{1,2}_\e.
\end{equation*}
We use the same resolution spaces and this simple inclusion in the perturbative analysis in section \ref{section3} in this paper.

Slightly later and independently, Bejenaru \cite{Be} also realized that the stereographic model can be used for perturbative analysis, and proved local well-posedness for small data in $H^\sigma$, in the full subcritical range $\sigma>d/2$, $d\geq 2$. In the stereographic model Bejenaru observed, apparently for the first time in the setting of Schr\"{o}dinger maps, that the gradient part of the nonlinearity has a certain null structure  (similar to the null structure of wave maps, observed by S. Klainerman).\footnote{This null structure was not observed in the earlier paper of Ionescu--Kenig \cite{IoKe2}; without this null structure the restriction $\sigma>(d+1)/2$ in \cite{IoKe2} is necessary for the perturbative argument.} The resolution spaces used in \cite{Be} for the perturbative argument are different from those of \cite{IoKe2}; these resolution spaces are based on the construction of suitably normalized wave packets, and had been previously used by Bejenaru in other subcritical problems (see \cite{Be2} and the references therein).

In \cite{IoKe3} Ionescu--Kenig proved the first global (in time) well-posedness theorem for small data in the critical Besov spaces $\dot{B}^{d/2}_Q$, in dimensions $d\geq 3$, using certain technical modifications of the resolution spaces of \cite{IoKe2} and the null structure observed in \cite{Be}. As explained in \cite{IoKe3}, the main difficulty in proving this result in dimension $d=2$ is the logarithmic failure of the scale-invariant $L^{2,\infty}_\e$ estimate. 

Unlike its Besov analogue, the condition $\|s_0-Q\|_{\dot{H}^{d/2}}\ll 1$ in Theorem \ref{Main1} does not guarantee that the data $s_0$ takes values in a small neighborhood of $Q$. Because of this, the stereographic model used in \cite{IoKe2}, \cite{Be}, and \cite{IoKe3} is not relevant, and it does not appear possible to prove Theorem \ref{Main1} using a direct perturbative construction. We construct the solution $s$ indirectly, using {\it{ a priori }} estimates: we start with a solution $s\in C([-T,T]:H^{\infty}_Q)$ of \eqref{Sch1}, where $T=T( \|s_0\|_{H^{\sigma_0}_Q})>0$, $\sigma_0$ sufficiently large, and transfer the quantitative bounds on the function $s$ at time $0$ to suitable quantitative bounds on the functions $\psi_m$ at time $0$ (the functions $\psi_m$ are solutions of the modified Schr\"{o}dinger map equations, see section \ref{gauge}). Then we study the modified Schr\"{o}dinger map equations perturbatively, and prove uniform quantitative bounds on the functions $\psi_m$ at all times $t\in[-T,T]$. Finally, we transfer these bounds back to the solution $s$; this gives uniform quantitative bounds on $s$ at all times $t\in[-T,T]$, which allow us to extend the solution $s$ up to time $T=1$. By scaling, we can construct a global solution.

The rest of the paper is organized as follows: in section \ref{gauge} we explain how to derive the modified Schr\"{o}dinger map equations (MSM)\footnote{The MSM were first derived in \cite{ChShUh}, using orthonormal frames, and \cite{NaStUh}, using the stereographic projection.}, and prove quantitative bounds on the solutions $\psi_m$ of the MSM at time $t=0$. In section \ref{section3} we use a perturbative argument and the resolution spaces defined in \cite{IoKe2} (and some of their properties) to prove bounds on the solutions $\psi_m$ of the MSM on the time interval $[-T,T]$. The proofs of some of the technical nonlinear bounds are deferred to section \ref{section5}. In section \ref{section4} we transfer the bounds on $\psi_m$ to a priori bounds on solution $s$ of \eqref{Sch1}, and use a local existence theorem to close the argument.

We will always assume in the rest of the paper that $d\geq 3$ (we have not constructed yet suitable resolution spaces in dimension $d=2$). In subsection \ref{section3.3} and sections \ref{section4} and \ref{section5} we assume the stronger restriction $d\geq 4$; the reason for this restriction is mostly technical, as it leads to simple proofs of the nonlinear estimates in Lemma \ref{Lemmaq5}. In many estimates, we will use the letter $C$ to denote constants that may depend only on the dimension $d$.

We would like to thank S. Klainerman, I. Rodnianski, J. Shatah, and T. Tao for several useful discussions.

\section{The modified Schr\"{o}dinger map}\label{gauge}

In this section we give a self-contained derivation of the modified Schr\"{o}dinger map equations, using orthonormal frames\footnote{This elementary construction was suggested to us by T. Tao.}. In the context of wave maps, orthonormal frames have been used in \cite{ChTa}, \cite{ShSt}, \cite{KlRo}, \cite{NaStUh3} etc. In the context of Schr\"{o}dinger maps, orthonormal frames (on the pullback of $T^\ast M$ under the solution $s$) have been used for the first time in \cite{ChShUh} to construct the modified Schr\"{o}dinger map equations. See also \cite{Ga}. Complete expositions of this construction have been presented by J. Shatah on several occasions.   

In this section we assume $d\geq 3$ (some technical changes are needed in dimension $d=2$, but we will not discuss them here).

\subsection{A topological construction}\label{section2.1}

Assume $n\in[1,\infty)\cap\Z$, $a_1,\ldots,a_n\in[0,\infty)$, and let
\begin{equation*}
\mathcal{D}^n=[-a_1,a_1]\times\ldots\times[-a_n,a_n].
\end{equation*}
For $n=0$ let $\mathcal{D}^0=\{0\}$. 

\newtheorem{Lemmag1}{Lemma}[section]
\begin{Lemmag1}\label{Lemmag1} Assume $n\geq 0$ and $s:\mathcal{D}^n\to\mathbb{S}^2$ is a continuous function. Then there is a continuous function $v:\mathcal{D}^n\to\mathbb{S}^2$ with the property that
\begin{equation*}
s(x)\cdot v(x)=0\text{ for any }x\in\mathcal{D}^n.
\end{equation*}
\end{Lemmag1}

\begin{proof}[Proof of Lemma \ref{Lemmag1}] We argue by induction over $n$ (the
case $n=0$ is trivial). Since $s$ is continuous, there is $\epsilon>0$ with the
property that
\begin{equation}\label{t1}
|s(x)-s(y)|\leq 2^{-10}\text{ for any }x,y\in\mathcal{D}^n\text{ with }|x-y|\leq\epsilon.
\end{equation}
For $x\in\mathcal{D}^n$ we write $x=(x',x_n)\in\mathcal{D}^{n-1}\times[-a_n,a_n]$. For any $b\in[-a_n,a_n]$ let $\mathcal{D}^n_b=\mathcal{D}^{n-1}\times[-a_n,b]=\{x=(x',x_n)\in\mathcal{D}^n:x_n\in[-a_n,b]\}$. By the induction hypothesis, we can define $v:\mathcal{D}^n_{-a_n}\to\mathbb{S}^2$ continuous such that
\begin{equation*}
s(x)\cdot v(x)=0\text{ for any }x\in\mathcal{D}^n_{-a_n}.
\end{equation*}
We extend now the function $v$ to $\mathcal{D}^n$. With $\epsilon$ as  in \eqref{t1}, it suffices to prove that if $b,b'\in[-a_n,a_n]$, $0\leq b'-b\leq\epsilon$, $v:\mathcal{D}^n_b\to\mathbb{S}^2$ is continuous, and $s(x)\cdot v(x)=0$ for any $x\in\mathcal{D}_b^n$, then $v$ can be extended to a continuous function $\widetilde{v}:\mathcal{D}^n_{b'}\to\mathbb{S}^2$ such that $s(x)\cdot \widetilde{v}(x)=0$ for any $x\in\mathcal{D}^n_{b'}$.

Let
\begin{equation}\label{t4}
\mathcal{R}=\{(u_1,u_2)\in\mathbb{R}^3\times\mathbb{R}^3:\,|u_1|,|u_2|\in(1/2,2)\text{ and }|u_1\cdot u_2|< 2^{-5}\},
\end{equation}
and let $N:\mathcal{R}\to\mathbb{S}^2$ denote the smooth function
\begin{equation}\label{t5}
N[u_1,u_2]=\frac{u_1-((u_1\cdot u_2)/|u_2|^2)\,u_2}{|u_1-((u_1\cdot u_2)/|u_2|^2)\,u_2|}.
\end{equation}
So $N[u_1,u_2]$ is a unit vector orthogonal to $u_2$ in the plane generated by the vectors $u_1$  and $u_2$. We construct now the extension $\widetilde{v}:\mathcal{D}^n_{b'}\to\mathbb{S}^2$. For $x'\in\mathcal{D}^{n-1}$ and $x_n\in[-a_n,b']$ let
\begin{equation*}
\widetilde{v}(x',x_n)=
\begin{cases}
N[v(x',b),s(x',x_n)]&\text{ if }x_n\in[b,b'];\\
v(x',x_n)&\text{ if }x_n\in[-a_n,b].
\end{cases}
\end{equation*}
In view of \eqref{t1}, the function $\widetilde{v}:\mathcal{D}^n_{b'}\to\mathbb{S}^2$ is well-defined, continuous, and $s(x)\cdot \widetilde{v}(x)=0$ for any $x\in\mathcal{D}_{b'}^n$. This completes the proof of Lemma \ref{Lemmag1}.
\end{proof}

\newtheorem{Lemmag2}[Lemmag1]{Lemma}
\begin{Lemmag2}\label{Lemmag2} Assume $T\in[0,2]$, $Q,Q'\in\mathbb{S}^2$, $Q\cdot Q'=0$, and $s:\mathbb{R}^d\times[-T,T]\to\mathbb{S}^2$ is a continuous function with the property that
\begin{equation*}
\lim_{x\to\infty}s(x,t)=Q\text{ uniformly in }t\in[-T,T].
\end{equation*}
Then there is a continuous function $v:\mathbb{R}^d\times[-T,T]\to\mathbb{S}^2$ with the property that
\begin{equation*}
\begin{cases}
&s(x,t)\cdot v(x,t)=0\text{ for any }(x,t)\in\mathbb{R}^d\times[-T,T];\\
&\lim\limits_{x\to\infty}v(x,t)=Q'\text{ uniformly in }t\in[-T,T].
\end{cases}
\end{equation*}
\end{Lemmag2}

\begin{proof}[Proof of Lemma \ref{Lemmag2}] We fix $R>0$ such that
\begin{equation*}
|s(x,t)-Q|\leq 2^{-10}\text{ if }|x|\geq R \text{ and }t\in[-T,T].
\end{equation*}
Using Lemma \ref{Lemmag1}, we can define a continuous function $v_0:B_R\times[-T,T]\to\mathbb{S}^2$ such that $s(x,t)\cdot v_0(x,t)=0$ for $(x,t)\in B_R\times[-T,T]$, where $B_R=\{x\in\mathbb{R}^d:|x|\leq R\}$. Let $S_R=\{x\in\mathbb{R}^d:|y|=R\}$ and $\mathbb{S}^1_Q=\{x\in\mathbb{S}^2:x\cdot Q=0\}$. We define the continuous function 
\begin{equation*}
w:S_R\times[-T,T]\to\mathbb{S}^1_Q,\,\,\,w(y,t)=\frac{(s(y,t)\cdot Q)v_0(y,t)-(v_0(y,t)\cdot Q)s(y,t)}{|(s(y,t)\cdot Q)v_0(y,t)-(v_0(y,t)\cdot Q)s(y,t)|},
\end{equation*} 
so $w(y,t)$ is a vector in $\mathbb{S}^1_Q$ and in the plane generated by $s(y,t)$ and $v_0(y,t)$. Since $d\geq 3$, the space $S_R\times[-T,T]$ is simply connected (and compact), thus the function $w$ is  homotopic to a constant function. Thus there is a continuous function 
\begin{equation*}
\widetilde{w}:S_R\times[-T,T]\times[1,2]\to\mathbb{S}^1_Q\,\text{ such that }\,\widetilde{w}(y,t,1)=w(y,t)\text{  and }\widetilde{w}(y,t,2)\equiv Q'.
\end{equation*}

With $N$ is as  in \eqref{t5}, we define
\begin{equation*}
v_1(x,t)=N[\widetilde{w}(Rx/|x|,t,|x|/R),s(x,t)]
\end{equation*}
for $|x|\in[R,2R]$, and
\begin{equation*}
v_2(x,t)=N[Q',s(x,t)]
\end{equation*}
for $|x|\geq 2R$. The function $v$ in Lemma \ref{Lemmag2} is obtained by gluing the functions $v_0$, $v_1$, and $v_2$.
\end{proof}

\subsection{Derivation of the modified Schr\"{o}dinger map equations}\label{section2.2}

Assume now that $T\in[0,1]$, $Q,Q'\in\mathbb{S}^2$, and $Q\cdot Q'=0$. Assume that 
\begin{equation}\label{t10}
\begin{cases}
&s\in C([-T,T]:H_Q^{\infty});\\
&\partial_ts\in C([-T,T]:H^{\infty}).
\end{cases}
\end{equation}
We extend the function $s$ to a function $\widetilde{s}\in C([-T-1,T+1]:H_Q^{\infty})$ by setting $\widetilde{s}(.,t)=s(.,T)$ if $t\in[T,T+1]$ and $\widetilde{s}(.,t)=s(.,-T)$ if $t\in[-T-1,-T]$. Clearly, the function $\widetilde{s}:\mathbb{R}^d\times[-T-1,T+1]\to\mathbb{S}^2$ is continuous and $\lim_{x\to\infty}\widetilde{s}(x,t)=Q$ uniformly in $t$. We apply Lemma \ref{Lemmag2} to construct a continuous function $\widetilde{v}:\mathbb{R}^d\times[-T-1,T+1]\to\mathbb{S}^2$ such that $\widetilde{s}\cdot\widetilde{v}\equiv 0$ and $\lim_{x\to\infty}\widetilde{v}(x,t)=Q'$ uniformly in $t$. 

We regularize now the function $\widetilde{v}$. Let $\varphi:\mathbb{R}^d\times\mathbb{R}\to[0,\infty)$ denote a smooth function supported in the ball $\{(x,t):|x|^2+t^2\leq 1\}$ with $\int_{\mathbb{R}^d\times\mathbb{R}}\varphi\,dxdt=1$. Since $\widetilde{v}$ is a uniformly continuous function, there is $\epsilon=\epsilon(\widetilde{v})$  with the property that
\begin{equation*}
|\widetilde{v}(x,t)-(\widetilde{v}\ast\varphi_\epsilon)(x,t)|\leq 2^{-20}\text{ for any }(x,t)\in\mathbb{R}^d\times[-T-1/2,T+1/2],
\end{equation*}
where $\varphi_\epsilon(x,t)=\epsilon^{-d-1}\varphi(x/ \epsilon,t/ \epsilon)$. Using a partition of $1$, we replace smoothly $(\widetilde{v}\ast\varphi_\epsilon)(x,t)$ with $Q'$ for $|x|$ large enough. Thus we have constructed a smooth function $v':\mathbb{R}^d\times(-T-1/2,T+1/2)\to\mathbb{R}^3$ with the properties
\begin{equation}\label{t11}
\begin{cases}
&|v'(x,t)|\in[1-2^{-10},1+2^{-10}]\text{ for any }(x,t)\in\mathbb{R}^d\times[-T,T];\\
&|v'(x,t)\cdot s(x,t)|\leq 2^{-10} \text{ for any }(x,t)\in\mathbb{R}^d\times[-T,T];\\
&v'(x,t)=Q'\text{ for }|x|\text{ large enough and }t\in[-T,T].
\end{cases}
\end{equation}

With $N$ as in \eqref{t5}, we define
\begin{equation*}
v(x,t)=N[v'(x,t),s(x,t)].
\end{equation*}
In view of \eqref{t11}, the continuous function $v:\mathbb{R}^d\times[-T,T]\to\mathbb{S}^2$ is well-defined, $s(x,t)\cdot v(x,t)\equiv 0$, and
\begin{equation}\label{t12}
\begin{cases}
&\partial_mv\in C([-T,T]:H^{\infty})\text{ for }m=1,\ldots,d;\\
&\partial_tv\in C([-T,T]:H^{\infty}).
\end{cases}
\end{equation}

Given $s$ as in \eqref{t10} and $v$ as in \eqref{t12}, we define
\begin{equation*}
w(x,t)=s(x,t)\times v(x,t).
\end{equation*}
Since $H^\sigma$ is an algebra for $\sigma>d/2$, we have
\begin{equation}\label{t13}
\begin{cases}
&\partial_mw\in C([-T,T]:H^{\infty})\text{ for }m=1,\ldots,d;\\
&\partial_tw\in C([-T,T]:H^{\infty }).
\end{cases}
\end{equation}
To summarize, given a function $s$ as in \eqref{t10} we have constructed continuous functions $v,w:\mathbb{R}^d\times[-T,T]\to\mathbb{S}^2$ such that $s\cdot v=s\cdot w=v\cdot w\equiv 0$, and \eqref{t12} and \eqref{t13} hold.

We use now the functions $v$ and $w$ to construct a suitable Coulomb gauge. Let
\begin{equation*}
A_m=(\partial_mv)\cdot w=-(\partial_mw)\cdot v\text{ for }m=1,\ldots,d.
\end{equation*}
Clearly, the functions $A_m$ are real-valued,
\begin{equation}\label{t20}
A_m\in C([-T,T]:H^{\infty})\text{ and }\partial_tA_m\in C([-T,T]:H^{\infty}). 
\end{equation}
We would like to modify the functions $v$ and $w$ such that $\sum_{m=1}^d\partial_mA_m\equiv 0$. Let
\begin{equation*}
\begin{cases}
&v'=(\cos \chi)v+(\sin\chi)w;\\
&w'=(-\sin\chi)v+(\cos\chi)w,
\end{cases}
\end{equation*}
for some function $\chi:\mathbb{R}^d\times[-T,T]\to\mathbb{R}$ to be determined. Then, using the orthonormality of $v$ and $w$ (which gives $\partial_mv\cdot v=\partial_mw\cdot w\equiv 0$),
\begin{equation*}
A'_m=(\partial_mv')\cdot w'=A_m+\partial_m\chi.
\end{equation*}
The condition $\sum_{m=1}^d\partial_mA'_m\equiv 0$ gives
\begin{equation*}
\Delta\chi=-\sum_{m=1}^d\partial_mA_m.
\end{equation*}
Thus we define $\chi$ by the formula
\begin{equation*}
\chi(x,t)=c\int_{\mathbb{R}^d}e^{ix\cdot \xi}|\xi|^{-2}\sum_{m=1}^d(i\xi_m)\,\mathcal{F}_{(d)}(A_m)(\xi,t)\,d\xi.
\end{equation*}
The integral defining the function $\chi$ converges absolutely since $A_m\in C([-T,T]:H^{\infty})$ and $d\geq 3$. Using \eqref{t20}, it follows that $\chi:\mathbb{R}^d\times[-T,T]\to\mathbb{R}$ is a bounded, continuous function, $\partial_m\chi\in C([-T,T]:H^{\infty})$ and $\partial_t\chi\in C([-T,T]:H^{\infty })$. To summarize, we proved the following proposition:

\newtheorem{Lemmag3}[Lemmag1]{Proposition}
\begin{Lemmag3}\label{Lemmag3} Assume $T\in[0,1]$, $Q\in\mathbb{S}^2$, and
\begin{equation}\label{t45}
\begin{cases}
&s\in C([-T,T]:H_Q^{\infty});\\
&\partial_ts\in C([-T,T]:H^{\infty}).
\end{cases}
\end{equation}
Then there are continuous functions $v,w:\mathbb{R}^d\times[-T,T]\to\mathbb{S}^2$, $s\cdot v\equiv 0$, $w=s\times v$, such that
\begin{equation}\label{t30}
\partial_mv,\partial_mw\in C([-T,T]:H^{\infty})\text{ for }m=0,1,\ldots,d,\\
\end{equation}
where $\partial_0=\partial_t$. In addition,
\begin{equation}\label{t40}
\text{ if }A_m=(\partial_mv)\cdot w\text{ for }m=1,\ldots,d,\text{ then }\sum_{j=1}^d\partial_mA_m\equiv 0.
\end{equation}
\end{Lemmag3}

Assume now that $s,v,w$ are as in Proposition \ref{Lemmag3}. In addition to the functions $A_m$, we define the continuous functions $\psi_m:\mathbb{R}^d\times[-T,T]\to\mathbb{C}$, $m=1,\ldots,d$,
\begin{equation}\label{t41}
\psi_m=(\partial_ms)\cdot v+i(\partial_ms)\cdot w.
\end{equation}
Let $\partial_0=\partial_t$. We also define the continuous functions $A_0:\mathbb{R}^d\times[-T,T]\to\mathbb{R}$ and $\psi_0:\mathbb{R}^d\times[-T,T]\to\mathbb{C}$,
\begin{equation}\label{t42}
\begin{cases}
&\psi_0=(\partial_0s)\cdot v+i(\partial_0s)\cdot w;\\
&A_0=(\partial_0v)\cdot w=-(\partial_0w)\cdot v.
\end{cases}
\end{equation}
Clearly, $\psi_m,A_m\in C([-T,T]:H^\infty)$ for $m=0,1,\ldots,d$, and $\partial_t\psi_m,\partial_tA_m\in C([-T,T]:H^\infty)$ for $m=1,\ldots,d$. In view of the orthonormality of $s,v,w$, for $m=0,1,\ldots,d$
\begin{equation}\label{t49}
\begin{cases}
&\partial_ms=\Re(\psi_m)v+\Im(\psi_m)w;\\
&\partial_mv=-\Re(\psi_m)s+A_mw;\\
&\partial_mw=-\Im(\psi_m)s-A_mv.
\end{cases}
\end{equation}

A direct computation using the orthonormality of $s,v,w$ gives
\begin{equation}\label{t43}
(\partial_l+iA_l)\psi_m=(\partial_m+iA_m)\psi_l\text{ for any }m,l=0,1,\ldots,d.
\end{equation}
A direct computation also shows that
\begin{equation}\label{t51}
\partial_lA_m-\partial_mA_l=\Im(\psi_l\,\overline{\psi}_m)\text{ for any }m,l=0,1,\ldots,d.
\end{equation}
We combine these identities with the Coulomb gauge condition $\sum_{m=1}^d\partial_mA_m\equiv 0$ and solve the div-curl system for each $t$ fixed. The result is
\begin{equation}\label{t50}
\Delta A_m=-\sum_{l=1}^d\partial_l[\Im(\psi_m\,\overline{\psi}_l)]\text{ for }m=1,\ldots,d.
\end{equation} 
Thus, using \eqref{t50}, for $m=1,\ldots,d$,
\begin{equation}\label{t52}
A_m=\nabla^{-1}\big[\sum_{l=1}^dR_l[\Im(\psi_m\,\overline{\psi}_l)]\big],
\end{equation}
where $R_l$ denotes the Riesz transform defined by the Fourier multiplier $\xi\to i\xi_l/|\xi|$ and $\nabla^{-1}$ is the operator defined by the Fourier multiplier $\xi\to|\xi|^{-1}$.

Assume now that the function $s$ satisfies the identity
\begin{equation}\label{t60}
\partial_ts=s\times\Delta s\,\text{ on }\,\mathbb{R}^d\times[-T,T],
\end{equation}
in addition to \eqref{t45}. For $m=0,1,\ldots,d$ we define the covariant derivatives $D_m=\partial_m+iA_m$. Using the definition,
\begin{equation*}
\psi_0=(s\times\Delta s)\cdot v+i(s\times\Delta s)\cdot w.
\end{equation*}
In addition, using \eqref{t49},
\begin{equation*}
\partial_m^2s=\big(\partial_m\Re(\psi_m)-A_m\cdot \Im(\psi_m)\big)v+\big(\partial_m\Im(\psi_m)+A_m\cdot \Re(\psi_m)\big)w-|\psi_m|^2s.
\end{equation*}
Thus, using $s\times v=w$, $s\times w=-v$,
\begin{equation}\label{t61}
\begin{split}
\psi_0&=-\sum_{m=1}^d\big(\partial_m\Im(\psi_m)+A_m\cdot \Re(\psi_m)\big)+i\sum_{m=1}^d\big(\partial_m\Re(\psi_m)-A_m\cdot \Im(\psi_m)\big)\\
&=i\sum_{m=1}^dD_m\psi_m.
\end{split}
\end{equation}

We use now \eqref{t43} and \eqref{t51} to convert \eqref{t61} into a nonlinear Schr\"{o}dinger equation. We rewrite the identities \eqref{t43} and \eqref{t51} in the form
\begin{equation*}
\begin{cases}
&D_l\psi_m=D_m\psi_l\,\text{ for any }m,l=0,1,\ldots,d;\\
&D_lD_mf-D_mD_lf=i\Im(\psi_l\overline{\psi}_m)f\,\text{ for any }m,l=0,1,\ldots,d.
\end{cases}
\end{equation*}
Thus, using \eqref{t61}, for $m=1,\ldots,d$,
\begin{equation*}
\begin{split}
D_0\psi_m&=D_m\psi_0=i\sum_{l=1}^dD_mD_l\psi_l=i\sum_{l=1}^dD_lD_m\psi_l-\sum_{l=1}^d\Im(\psi_m\overline{\psi}_l)\psi_l\\
&=i\sum_{l=1}^dD_lD_l\psi_m-\sum_{l=1}^d\Im(\psi_m\overline{\psi}_l)\psi_l.
\end{split}
\end{equation*}
Thus, using again \eqref{t40}, for $m=1,\ldots,d$,
\begin{equation}\label{t62}
(i\partial_t+\Delta_x)\psi_m=-2i\sum_{l=1}^dA_l\cdot \partial_l\psi_m+\big(A_0+\sum_{l=1}^dA_l^2\big)\psi_m-i\sum_{l=1}^d\Im(\psi_m\overline{\psi}_l)\psi_l.
\end{equation}

We find now the coefficient $A_0$. Using \eqref{t51} and \eqref{t40},
\begin{equation}\label{t70}
\Delta A_0=\sum_{l=1}^d\partial_l(\partial_0A_l+\Im(\psi_l\overline{\psi}_0))=\sum_{l=1}^d\partial_l\,\Im(\psi_l\overline{\psi}_0).
\end{equation}
Using \eqref{t61}, \eqref{t43} and the identity $\overline{\psi}_l\cdot D_m\psi_m=\partial_m(\overline{\psi}_l\psi_m)-\psi_m\cdot\overline{D_m\psi_l}$,
\begin{equation*}
\begin{split}
\Im(\psi_l\overline{\psi}_0)&=-\sum_{m=1}^d\Re(\overline{\psi}_l\cdot D_m\psi_m)=-\sum_{m=1}^d\partial_m\Re(\overline{\psi}_l\psi_m)+\sum_{m=1}^d\Re(\psi_m\cdot\overline{D_m\psi_l})\\
&=-\sum_{m=1}^d\partial_m\Re(\overline{\psi}_l\psi_m)+\frac{1}{2}\partial_l\big(\sum_{m=1}^d\psi_m\overline{\psi}_m\big).
\end{split}
\end{equation*}
It  follows from \eqref{t70} that
\begin{equation*}
\Delta A_0=-\sum_{m,l=1}^d\partial_l\partial_m\Re(\overline{\psi}_l\psi_m)+\frac{1}{2}\Delta\big(\sum_{m=1}^d\psi_m\overline{\psi}_m\big).
\end{equation*}
Thus
\begin{equation}\label{t75}
A_0=\sum_{m,l=1}^dR_lR_m\big(\Re(\overline{\psi}_l\psi_m)\big)+\frac{1}{2}\sum_{m=1}^d\psi_m\overline{\psi}_m.
\end{equation}

\newtheorem{Lemmag4}[Lemmag1]{Proposition}
\begin{Lemmag4}\label{Lemmag4} Assume $s,v,w$, and $A_m$, $m=1,\ldots,d$ are as in Proposition \ref{Lemmag3}. Assume in addition that the function $s$ satisfies the identity
\begin{equation*}
\partial_ts=s\times\Delta s\,\text{ on }\,\mathbb{R}^d\times[-T,T].
\end{equation*}
For $m=1,\ldots,d$ let
\begin{equation}\label{t80}
\psi_m=(\partial_ms)\cdot v+i(\partial_ms)\cdot w\text{ on }\mathbb{R}^d\times[-T,T].
\end{equation}
Then $\psi_m,A_m,\partial_t\psi_m,\partial_tA_m\in C([-T,T]:H^{\infty})$ and 
\begin{equation}\label{t81}
\begin{cases}
&(\partial_l+iA_l)\psi_m=(\partial_m+iA_m)\psi_l\text{ for any }m,l=1,\ldots,d;\\
&A_m=\nabla^{-1}\big[\sum_{l=1}^dR_l[\Im(\psi_m\,\overline{\psi}_l)]\big]\text{ for any }m=1,\ldots,d,
\end{cases}
\end{equation}
where $R_l$ denotes the Riesz transform defined by the Fourier multiplier $\xi\to i\xi_l/|\xi|$ and $\nabla^{-1}$ is the operator defined by the Fourier multiplier $\xi\to |\xi|^{-1}$. In addition, the functions $\psi_m$ satisfy the system of nonlinear Schr\"{o}dinger equations
\begin{equation}\label{t82}
(i\partial_t+\Delta_x)\psi_m=-2i\sum_{l=1}^dA_l\cdot \partial_l\psi_m+\big(A_0+\sum_{l=1}^dA_l^2\big)\psi_m+i\sum_{l=1}^d\Im(\psi_l\overline{\psi}_m)\psi_l,
\end{equation}
for $m=1,\ldots,d$, where
\begin{equation}\label{t83}
A_0=\sum_{l,l'=1}^dR_lR_{l'}\big(\Re(\overline{\psi}_l\psi_{l'})\big)+\frac{1}{2}\sum_{l=1}^d\psi_l\overline{\psi}_l.
\end{equation}
\end{Lemmag4}

\subsection{A quantitative estimate}\label{section2.3}
We prove now quantitative estimates for the functions $\psi_m$.

\newtheorem{Lemmag5}[Lemmag1]{Lemma}
\begin{Lemmag5}\label{Lemmag5} With the notation in Propositions \ref{Lemmag3} and \ref{Lemmag4}, if the function $s_0(x)=s(x,0)$ has the additional property $\|s_0-Q\|_{\dot{H}^{d/2}}\leq1$ and $\sigma_0=d+10$ then, for $m=1,\ldots,d$,
\begin{equation}\label{t90}
\begin{cases}
&\|\psi_m(.,0)\|_{\dot{H}^{(d-2)/2}}\leq C\cdot \|s_0-Q\|_{\dot{H}^{d/2}};\\
&\|\psi_m(.,0)\|_{H^{\sigma'-1}}\leq C( \|s_0\|_{H^{\sigma'}_Q})\text{ for any }\sigma'\in[1,\sigma_0]\cap\Z.
\end{cases}
\end{equation}
\end{Lemmag5}

\begin{proof}[Proof of Lemma \ref{Lemmag5}] The main difficulty is that our construction does not give effective control of the Sobolev norms of $v$ and $w$ in terms of the norms of $s$. We argue indirectly, using a bootstrap argument and the identities \eqref{t49}, \eqref{t80}, and \eqref{t81}. For $\sigma\in [-1,\infty)$ let $\nabla^\sigma$ denote the operator (acting on functions in $H^{\infty}$) defined by the Fourier multiplier $\xi\to|\xi|^\sigma$. For $\sigma\in[-1/2,d/2]$ let $p_\sigma=d/(\sigma+1)$. Then, in view of the Sobolev imbedding theorem (recall $d\geq 3$),
\begin{equation}\label{t91}
\|\nabla^\sigma f\|_{L^{p_\sigma}}\leq C\|\nabla^{\sigma'}f\|_{L^{p_{\sigma'}}}\text{ if }-1/2\leq \sigma\leq\sigma'\leq d/2\text{ and }f\in H^{\infty}.
\end{equation}

Let $s_0(x)=s(x,0)$, $v_0(x)=v(x,0)$, $w_0(x)=w(x,0)$, $\psi_{m,0}(x)=\psi_m(x,0)$, and $A_{m,0}(x)=A_m(x,0)$, and let $\epsilon_0=\|s_0-Q\|_{\dot{H}^{d/2}}\leq 1$. To start our bootstrap argument, we use \eqref{t80}, \eqref{t91} and the fact that $|v_0|=|w_0|=1$ to obtain
\begin{equation*}
\|\psi_{m,0}\|_{L^{p_0}}\leq C\epsilon_0\text{ for }m=1,\ldots,d.
\end{equation*}
Then, using \eqref{t81}, 
\begin{equation*}
\|\nabla^1A_{m,0}\|_{L^{p_1}}\leq C\epsilon_0\text{ for }m=1,\ldots,d.
\end{equation*}
Thus, using \eqref{t91}, $\|A_{m,0}\|_{L^{p_0}}\leq C\epsilon_0$ for $m=1,\ldots,d$. We use now the identity \eqref{t49} and the fact that for $f\in H^{\infty }$
\begin{equation}\label{t92}
\|\nabla^nf\|_{L^p}\approx \sum_{n_1+\ldots+n_d=n}\|\partial_1^{n_1}\ldots\partial_d^{n_d}f\|_{L^p}\text{ if }n\in\mathbb{Z}_+\text{ and }p\in[p_{d/2},p_{-1/2}].
\end{equation}
Thus
\begin{equation*}
\|\nabla^1v_0\|_{L^{p_0}}+\|\nabla^1w_0\|_{L^{p_0}}\leq C\epsilon_0.
\end{equation*}
Therefore
\begin{equation}\label{t94}
\sum_{m=1}^d\|\psi_{m,0}\|_{L^{p_0}}+\sum_{m=1}^d\|\nabla^1A_{m,0}\|_{L^{p_1}}+\|\nabla^1v_0\|_{L^{p_0}}+\|\nabla^1w_0\|_{L^{p_0}}\leq C\epsilon_0.
\end{equation}

We prove now that
\begin{equation}\label{t96}
\sum_{m=1}^d\|\nabla^n\psi_{m,0}\|_{L^{p_n}}+\sum_{m=1}^d\|\nabla^{n+1}A_{m,0}\|_{L^{p_{n+1}}}+\|\nabla^{n+1}v_0\|_{L^{p_n}}+\|\nabla^{n+1}w_0\|_{L^{p_n}}\leq C\epsilon_0,
\end{equation}
for any $n\in\mathbb{Z}\cap[0,(d-2)/2]$. We argue by induction over $n$. The case $n=0$ was already proved in \eqref{t94}. Assume $n\geq 1$ and \eqref{t96} holds for any $n'\in[0,n-1]\cap\mathbb{Z}$. Using \eqref{t80}, \eqref{t92}, and the induction hypothesis
\begin{equation*}
\begin{split}
\|\nabla^n\psi_{m,0}\|_{L^{p_n}}&\leq C\|\nabla^{n+1}s_0\|_{L^{p_n}}\cdot \|v_0\|_{L^\infty}\\
&+C\sum_{n'=0}^{n-1}\|\nabla^{n-n'}s_0\|_{L^{p_{n-n'-1}}}\cdot \|\nabla^{n'+1}v_0\|_{L^{p_{n'}}},
\end{split}
\end{equation*}
which suffices to control the first term in the left-hand side of \eqref{t96}. For the second term, using \eqref{t81} and \eqref{t92},
\begin{equation*}
\|\nabla^{n+1}A_{m,0}\|_{L^{p_{n+1}}}\leq C\sum_{l,l'=1}^d\sum_{n'=0}^{n}\|\nabla^{n'}\psi_{l,0}\|_{L^{p_{n'}}}\cdot \|\nabla^{n-n'}\psi_{l',0}\|_{L^{p_{n-n'}}},
\end{equation*}
which suffices in view of the induction hypothesis and the bound on the first term proved before. The bound on the last two terms in the left-hand side of \eqref{t96} follows in a similar way, using \eqref{t49}, \eqref{t92}, and the bound on the first two terms.

If $d$ is even then \eqref{t96} suffices to prove the first inequality in \eqref{t90}, simply by taking $n=(d-2)/2$. If $d$ is odd, the bounds \eqref{t96} with $n=(d-3)/2$ and \eqref{t91} give
\begin{equation}\label{t98}
\|\nabla^{\sigma+1}v_0\|_{L^{p_\sigma}}+\|\nabla^{\sigma+1}w_0\|_{L^{p_\sigma}}\leq C\epsilon_0\text{ for }\sigma\in[-1/2,(d-3)/2].
\end{equation}
In view of the hypothesis and \eqref{t91}, we also have the bound
\begin{equation}\label{t99}
\|\nabla^{\sigma+1}s_0\|_{L^{p_\sigma}}\leq C\epsilon_0\text{ for }\sigma\in[-1/2,(d-2)/2].
\end{equation}
We need the following Leibniz rule (a particular case of \cite[Theorem A.8]{KePoVe}):
\begin{equation}\label{Le}
\|\nabla^{1/2}(fg)-g\nabla^{1/2}f\|_{L^2}\leq C\|\nabla^{1/2}g\|_{L^{q_1}}\cdot \|f\|_{L^{q_2}}
\end{equation}
if $1/q_1+1/q_2=1/2$ and $q_1,q_2\in[p_{d/2},p_{-1/2}]$. Then, using \eqref{t80} and \eqref{t92}
\begin{equation*}
\|\nabla^{(d-2)/2}\psi_{m,0}\|_{L^2}\leq C\sum_{u_0\in\{v_0,w_0\}}\sum_{n=0}^{(d-3)/2}\|\nabla^{1/2}(\partial_mD^ns_0\cdot D^{(d-3)/2-n}u_0)\|_{L^2},
\end{equation*}
where $D^n$ denotes any derivative of the form $\partial_1^{n_1}\ldots\partial_d^{n_d}$, with $n_1+\ldots+n_d=n$. The first inequality in \eqref{t90} then follows from \eqref{t98}, \eqref{t99}, \eqref{Le} and the fact that $|u_0|\equiv 1$.

For the second inequality in \eqref{t90}, we notice first that $\|\psi_{m,0}\|_{H^0}\leq C\cdot \|s_0\|_{H^{1}_Q}$, since $|v_0|=|w_0|\equiv 1$. In view of the first inequality in \eqref{t90}, we may assume $\sigma'\geq (d+1)/2$. We use a similar argument as before: the bootstrap inequality that replaces \eqref{t96} is
\begin{equation}\label{t96.1}
\begin{split}
\sum_{m=1}^d\|\nabla^n&\psi_{m,0}\|_{L^2\cap L^{p_{n-\sigma'+d/2}}}+\sum_{m=1}^d\|\nabla^{n}A_{m,0}\|_{L^2\cap L^{p_{n-\sigma'+d/2}}}\\
&+\sum_{u_0\in \{v_0,w_0\}}\|\nabla^{n+1}u_0\|_{L^2\cap L^{p_{n-\sigma'+d/2}}}\leq C( \|s_0\|_{H^{\sigma'}_Q}),
\end{split}
\end{equation}
for any $n\in[0,\sigma'-1]\cap\Z$,  where $p_\sigma=p_{-1/2}=2d$  if $\sigma\leq-1/2$. As before, the bound \eqref{t96.1} follows by induction over $n$, using the identities \eqref{t49}, \eqref{t80}, and \eqref{t81}, and the inequalities \eqref{t91}, \eqref{t92}, and
\begin{equation*}
\sum_{n_1+\ldots+n_d\leq \sigma'-(d+1)/2}||\partial_1^{n_1}\ldots\partial_d^{n_d}s_0||_{L^\infty}\leq C( \|s_0\|_{H^{\sigma'}_Q}).
\end{equation*}
The second inequality in \eqref{t90} follows from the bound \eqref{t96.1} with $n=\sigma'-1$. 
\end{proof}

\section{Perturbative analysis of the modified Schr\"{o}dinger map}\label{section3}

In this section we analyze the Schr\"{o}dinger map system derived in Propositions \ref{Lemmag3} and \ref{Lemmag4}. In the rest of this section we assume $d\geq 3$; this restriction is used implicitly in many estimates. 

\subsection{The resolution spaces and their properties}\label{section3.1}

In this subsection we define our main normed spaces and summarize some of their basic properties. These resolution spaces have been used in \cite{IoKe2} and, with slight modifications, in \cite{IoKe3}, and we will refer to these papers for most of the proofs. 

Let $\mathcal{F}$ and $\mathcal{F}^{-1}$ denote the Fourier transform and the inverse Fourier transform operators on $L^2(\mathbb{R}^{d+1})$. For $l=1,\ldots,d$ let $\mathcal{F}_{(l)}$ and $\mathcal{F}_{(l)}^{-1}$ denote the Fourier transform and the inverse Fourier transform operators on $L^2(\mathbb{R}^l)$. We fix $\eta_0:\mathbb{R}\to[0,1]$ a smooth even function supported in the set $\{\mu\in\mathbb{R}:|\mu|\leq 8/5\}$ and equal to $1$ in the set $\{\mu\in\mathbb{R}:|\mu|\leq 5/4\}$. Then we define $\eta_j:\mathbb{R}\to[0,1]$, $j=1,2,\ldots$,
\begin{equation}\label{gu1}
\eta_j(\mu)=\eta_0(\mu/2^j)-\eta_0(\mu/2^{j-1}),
\end{equation}
and $\eta_k^{(d)}:\mathbb{R}^d\to[0,1]$, $k\in\Z$,
\begin{equation}\label{gu1.1}
\eta_k^{(d)}(\xi )=\eta_0( |\xi|/2^k)-\eta_0( |\xi|/2^{k-1}).
\end{equation}
For $j\in\Z_+$, we also define $\eta_{\leq j}=\eta_0+\ldots+\eta_j$.

For $k\in\mathbb{Z}$ let $I_k^{(d)}=\{\xi\in\R^d:|\xi|\in[2^{k-1},2^{k+1}]\}$; for $j\in\Z_+$ let $I_j=\{\mu\in\R:|\mu|\in[2^{j-1},2^{j+1}]\}$ if $j\geq 1$ and $I_j=[-2,2]$ if $j=0$. For $k\in\Z$ and $j\in\Z_+$ let
\begin{equation*}
D_{k,j}=\{(\xi,\tau)\in\mathbb{R}^d\times\mathbb{R}:\xi \in I_k^{(d)}\text{ and }|\tau+|\xi|^2|\in I_j\}\text{ and }D_{k,\leq j}=\bigcup\limits_{0\leq j'\leq j}D_{k,j'}.
\end{equation*}

For $k\in\Z$ we define first the normed spaces
\begin{equation}\label{v1}
\begin{split}
X_k=\{f\in L^2(\mathbb{R}^d\times&\mathbb{R}):f\text{ supported in }I_k^{(d)}\times\mathbb{R}\text { and } \\
&\|f\|_{X_k}=\sum_{j=0}^\infty 2^{j/2}\|\eta_j(\tau+|\xi|^2)\cdot f\|_{L^2}<\infty\}.
\end{split}
\end{equation}
The spaces $X_k$ are not sufficient for our estimates, due to various logarithmic divergences. For any vector $\mathbf{e}\in\mathbb{S}^{d-1}$ let $$P_{\mathbf{e}}=\{\xi\in\mathbb{R}^d:\xi\cdot\mathbf{e}=0\}$$ with the induced Euclidean measure. For $p,q\in[1,\infty]$ we define the normed spaces $L^{p,q}_{\mathbf{e}}=L^{p,q}_{\mathbf{e}}(\mathbb{R}^d\times\mathbb{R})$,
\begin{equation}\label{vv1}
\begin{split}
L^{p,q}_{\mathbf{e}}&=\{f\in L^2(\mathbb{R}^d\times \mathbb{R}):\\
&\|f\|_{L^{p,q}_{\mathbf{e}}}=\Big[\int_{\mathbb{R}}\Big[\int_{P_\mathbf{e}\times \mathbb{R}}|f(r\mathbf{e}+v,t)|^q\,dvdt\Big]^{p/q}\,dr\Big]^{1/p}<\infty\}.
\end{split}
\end{equation}
For $k\in\Z$ and $j\in \mathbb{Z}_+$ let
\begin{equation*}
D_{k,j}^{\mathbf{e}}=\{(\xi,\tau)\in D_{k,j}:\xi\cdot\mathbf{e}\geq 2^{k-20}\}\text{ and }D_{k,\leq j}^{\mathbf{e}}=\bigcup_{0\leq j'\leq j}D_{k,j}^{\mathbf{e}}.
\end{equation*}
For $k\geq 100$ and $\mathbf{e}\in\mathbb{S}^{d-1}$, we define the normed spaces
\begin{equation}\label{v2}
\begin{split}
Y_k^{\mathbf{e}}&=\{f\in L^2(\mathbb{R}^d\times\mathbb{R}):f\text{ supported in }D_{k,\leq 2k+10}^{\mathbf{e}}\text { and } \\
&\|f\|_{Y_k^{\mathbf{e}}}=2^{-k/2}\|\mathcal{F}^{-1}[(\tau+|\xi|^2+i)\cdot f]\|_{L^{1,2}_{\mathbf{e}}}<\infty\}.
\end{split}
\end{equation}
For simplicity of notation, we also define $Y_k^{\mathbf{e}}=\{0\}$ for $k\leq 99$.

We fix $L=L(d)$ large and $\mathbf{e}_1,\ldots,\mathbf{e}_L\in\mathbb{S}^{d-1}$, $\mathbf{e}_l\neq \mathbf{e}_{l'}$ if $l\neq l'$, such that
\begin{equation}\label{vm4}
\text{ for any }\mathbf{e}\in\mathbb{S}^{d-1}\text{ there is }l\in\{1,\ldots,L\}\text{ such that }|\mathbf{e}-\mathbf{e}_l|\leq 2^{-100}.
\end{equation}
We assume in addition that if $\mathbf{e}\in\{\mathbf{e}_1,\ldots,\mathbf{e}_L\}$ then $-\mathbf{e}\in\{\mathbf{e}_1,\ldots,\mathbf{e}_L\}$. For $k\in\mathbb{Z}$ we define the normed spaces
\begin{equation}\label{v3'}
Z_k=X_k+Y_{k}^{\mathbf{e}_1}+\ldots+Y_k^{\mathbf{e}_L}.
\end{equation}
The spaces $Z_k$ are our main normed spaces.

For $k\in \Z_+$ let $\Xi_k=2^k\cdot\Z^d$. Let $\chi^{(1)}:\mathbb{R}\to[0,1]$ denote an even smooth function supported in the interval $[-2/3,2/3]$ with the property that
\begin{equation*}
\sum_{n\in\Z}\chi^{(1)}(\xi-n)\equiv 1\text{ on }\R.
\end{equation*}
Let $\chi:\mathbb{R}^d\to[0,1]$, $\chi(\xi)=\chi^{(1)}(\xi_1)\cdot\ldots\cdot \chi^{(1)}(\xi_d)$. For $k\in \Z_+$ and $n\in\Xi_k$ let
\begin{equation*}
\chi_{k,n}(\xi)=\chi((\xi-n)/2^k).
\end{equation*}
Clearly, $\sum_{n\in\Xi_k}\chi_{k,n}\equiv 1$ on $\R^d$.

We summarize now some of the main properties of the spaces $Z_k$.

\newtheorem{Lemmas1}{Proposition}[section]
\begin{Lemmas1}\label{Lemmas1}
(a) If $k\in\Z$, $m\in L^\infty(\mathbb{R}^d)$, $\mathcal{F}_{(d)}^{-1}(m)\in L^1(\mathbb{R}^d)$, and $f\in Z_k$,  then $m(\xi)\cdot f\in Z_k$ and
\begin{equation}\label{im1}
||m(\xi)\cdot f||_{Z_k}\leq C||\mathcal{F}_{(d)}^{-1}(m)||_{L^1(\mathbb{R}^d)}\cdot ||f||_{Z_k}.
\end{equation}

(b) If $k\in\Z$, $j\in\mathbb{Z}_+$ and $f\in Z_k$ then
\begin{equation}\label{im2}
\|f\cdot \eta_j(\tau+|\xi|^2)\|_{X_k}\leq C\|f\|_{Z_k}.
\end{equation}

(c) If $k\in\Z$, $j\in\Z_+$, and $f\in Z_k$ then
\begin{equation}\label{im3}
||\eta_{\leq j}(\tau+|\xi|^2)\cdot f||_{Z_k}\leq C||f||_{Z_k}.
\end{equation}

(d) If $k\in\Z$ and $f$ is supported in $D^{\e}_{k,\leq \infty}$ for some $\e\in\{\e_1,\ldots,\e_L\}$ then
\begin{equation}\label{im4}
||f||_{Z_k}\leq C 2^{-k/2}||\mathcal{F}^{-1}
[(\tau+|\xi|^2+i)\cdot f]||_{L^{1,2}_\e}.
\end{equation}

(e) (Energy estimate) If $k\in\Z$ and $f\in Z_k$ then
\begin{equation}\label{im5}
\sup_{t\in\mathbb{R}}\|\mathcal{F}^{-1}(f)(.,t)\|_{L^2_x}
\leq C\|f\|_{Z_k}.
\end{equation}

(f) (Localized maximal function estimate) If $k\in\Z$, $k'\in(-\infty,k+10d]\cap\Z$, $f\in Z_k$, and $\mathbf{e}'\in\mathbb{S}^{d-1}$ then
\begin{equation}\label{im7}
\big[\sum_{n\in\Xi_{k'}}||\mathcal{F}^{-1}(\chi_{k',n}(\xi)\cdot \widetilde{f})||_{L^{2,\infty}_{\e'}}^2\big]^{1/2}\leq C2^{(d-1)k/2}\cdot 2^{-(d-2)(k-k')/2}(1+|k-k'| )\cdot \|f\|_{Z_k},
\end{equation}
where $\mathcal{F}^{-1}(\widetilde{f})\in\{\mathcal{F}^{-1}(f),\overline{\mathcal{F}^{-1}(f)}\}$
.

(g) (Local smoothing estimate) If $k\in \Z$, $\mathbf{e}'\in\mathbb{S}^{d-1}$,
$l\in[-1,40]\cap\Z$, and $f\in Z_k$ then
\begin{equation}\label{im8}
\|\mathcal{F}^{-1}[\widetilde{f}\cdot
\eta_{1}(\xi\cdot\mathbf{e}'/2^{k-l})]\|_{L^{\infty,2}_{\mathbf{e}'}}\leq
C2^{-k/2}\|f\|_{Z_k},
\end{equation}
where $\mathcal{F}^{-1}(\widetilde{f})\in\{\mathcal{F}^{-1}(f),\overline{\mathcal{F}^{-1}(f)}\}$.
\end{Lemmas1}

The bound \eqref{im1} follows directly from the definitions. The bound \eqref{im2} is proved in \cite[Lemma 2.1]{IoKe3}. The bound \eqref{im3} is  proved in \cite[Lemma 2.3]{IoKe3}. The bound \eqref{im4} follows from the estimate (2.15) in \cite{IoKe3}. The energy estimate \eqref{im5} is proved in \cite[Lemma 2.2]{IoKe3}. The localized maximal function estimate \eqref{im7} follows from \cite[Lemma 4.1]{IoKe3} and \eqref{im2}. Finally, the local smoothing estimate \eqref{im8} is proved in \cite[Lemma 4.2]{IoKe3}.

The estimate in part (f) with $k'=k$ will often be referred to as the ``global \eqref{im7}''. For $k'\leq k-C$ we refer to this estimate as the ``localized \eqref{im7}''.

\subsection{Linear estimates}\label{section3.2} 
We fix a large constant $\sigma_0$, say
\begin{equation}\label{me1}
\sigma_0=d+10.
\end{equation}
For $\sigma\in[(d-2)/2,\sigma_0-1]$ we define the normed space
\begin{equation}\label{no5}
\begin{split}
\dot{F}^\sigma&=\{u\in C(\mathbb{R}:H^{\infty}):\\
&\|u\|_{\dot{F}^\sigma}=\big[\sum_{k\in\Z}(2^{2\sigma k}+2^{(d-2)k})\,\|\eta_k^{(d)}(\xi)\cdot \mathcal{F}(u)\|^2_{Z_k}\big]^{1/2}<\infty\}.
\end{split}
\end{equation}
For $\sigma\in[(d-2)/2,\sigma_0-1]$, $T\in[0,1]$, $u\in C([-T,T]:H^{\infty})$, and $T'\in[0,T]$ we define
\begin{equation}\label{no6.1}
E_{T'}(u)(t)=
\begin{cases}
u(t)&\text{ if }|t|\leq T';\\
0&\text{ if }|t|>T',
\end{cases}
\end{equation} 
and
\begin{equation}\label{no6}
\begin{split}
\|u\|_{\dot{N}^\sigma[-T',T']}=\big[\sum_{k\in\Z}(2^{2\sigma k}+2^{(d-2)k})\,\|\eta_k^{(d)}(\xi)\cdot (\tau+|\xi|^2+i)^{-1}\cdot \mathcal{F}(E_{T'}u)\|_{Z_k}^2\big]^{1/2}.
\end{split}
\end{equation}
The definition \eqref{v1} shows that if $k\in\Z$ and $f$ is supported in $I_k^{(d)}\times\R$ then
\begin{equation*}
\|(\tau+|\xi|^2+i)^{-1}\cdot f\|_{Z_k}\leq C\|f\|_{L^2},
\end{equation*}
thus, for $\sigma\in[(d-2)/2,\sigma_0-1]$ and $T_1,T_2\in[0,T]$
\begin{equation}\label{no6.2}
\big|\|u\|_{\dot{N}^\sigma[-T_1,T_1]}-\|u\|_{\dot{N}^\sigma[-T_2,T_2]}\big|\leq C|T_1-T_2|^{1/2}\cdot \sup_{t\in[-T,T]}\|u(.,t)\|_{H^\sigma}.
\end{equation}

For $\phi\in H^{\sigma}$ let $W(t)(\phi)\in
C(\mathbb{R}:H^{\sigma})$ denote the solution of the free Schr\"{o}dinger evolution. 

\newtheorem{Lemmaqq1}[Lemmas1]{Proposition}
\begin{Lemmaqq1}\label{Lemmaqq1}
If $\sigma\in [(d-2)/2,\sigma_0-1]$ and $\phi\in H^{\infty}$ then
\begin{equation*}
\|\eta_0(t)\cdot W(t)(\phi)\|_{\dot{F}^{\sigma}}\leq C(\|\phi\|_{\dot{H}^\sigma}+\|\phi\|_{\dot{H}^{(d-2)/2}}).
\end{equation*}
\end{Lemmaqq1}

See \cite[Lemma 3.1]{IoKe3} for the proof.

\newtheorem{Lemmaq3}[Lemmas1]{Proposition}
\begin{Lemmaq3}\label{Lemmaq3}
If $\sigma\in[(d-2)/2,\sigma_0-1]$, $T\in[0,1]$, and $u\in C([-T,T]:H^{\infty})$ then
\begin{equation*}
\Big|\Big|\eta_0(t)\cdot \int_0^tW(t-s)(E_T(u)(s))\,ds\Big|\Big|_{\dot{F}^{\sigma}}\leq C||u||_{\dot{N}^{\sigma}[-T,T]},
\end{equation*}
where $E_T(u)$ is defined in \eqref{no6.1}.
\end{Lemmaq3}

See \cite[Lemma 3.2]{IoKe3} for the proof.

\subsection{Nonlinear estimates}\label{section3.3} In this subsection we assume that $d\geq 4$. Assume that $T\in[0,1]$ and $\psi_m\in C([-T,T]:H^{\infty })$, $m=1,\ldots,d$. Let $\Psi=(\psi_1,\ldots,\psi_d)$ and define
\begin{equation}\label{bh1}
\begin{cases}
&A_0=\sum_{l,l'=1}^dR_lR_{l'}\big(\Re(\overline{\psi}_l\psi_{l'})\big)+\frac{1}{2}\sum_{l=1}^d\psi_l\overline{\psi}_l;\\
&A_m=\nabla^{-1}\big[\sum_{l=1}^dR_l[\Im(\psi_m\,\overline{\psi}_l)]\big]\text{ for any }m=1,\ldots,d,
\end{cases}
\end{equation}
and
\begin{equation}\label{bh2}
\mathcal{N}_m(\Psi)=-2i\sum_{l=1}^dA_l\cdot \partial_l\psi_m+\big(A_0+\sum_{l=1}^dA_l^2\big)\psi_m+i\sum_{l=1}^d\Im(\psi_l\overline{\psi}_m)\psi_l.
\end{equation}
Clearly, $A_m,\mathcal{N}_m(\Psi)\in C([-T,T]:H^{\infty})$ (recall that $d\geq 3$). We assume also that on $\mathbb{R}^d\times[-T,T]$ we have the integral equation
\begin{equation}\label{inte}
\psi_m(t)=W(t)(\psi_{m,0})+\int_0^tW(t-s)(\mathcal{N}_m(\Psi)(s))\,ds,
\end{equation}
where $\psi_{m,0}=\psi_m(0)$. In dimensions $d\geq 4$ we will not need the compatibility conditions
\begin{equation*}
(\partial_l+iA_l)\psi_m=(\partial_m+iA_m)\psi_l\text{ for any }m,l=1,\ldots,d.
\end{equation*}

We define the extensions $\widetilde{E}_T(\psi_m)\in C(\R:H^{\infty})$, $m=1,\ldots,d$,
\begin{equation}\label{ty4}
\widetilde{E}_T(\psi_m)(t)=\eta_0(t)\cdot W(t)(\psi_{m,0})+\eta_0(t)\cdot \int_0^tW(t-s)(E_T(\mathcal{N}_m(\Psi))(s))\,ds.
\end{equation}
Using Propositions \ref{Lemmaqq1} and \ref{Lemmaq3}, for $\sigma\in[(d-2)/2,\sigma_0-1]$
\begin{equation*}
\|\widetilde{E}_T(\psi_m)\|_{\dot{F}^\sigma}\leq C\cdot ( \|\psi_{m,0}\|_{\dot{H}^\sigma\cap\dot{H}^{(d-2)/2}}+\|\mathcal{N}_m(\Psi)\|_{\dot{N}^\sigma[-T,T]}).
\end{equation*}
Let $\widetilde{E}_T(\Psi)=(\widetilde{E}_T(\psi_1),\ldots,\widetilde{E}_T(\psi_d))$. For $\sigma\in[(d-2)/2,\sigma_0-1]$ let
\begin{equation}\label{bh2.1}
\|\widetilde{E}_T(\Psi)\|_{\dot{F}^\sigma}=\sum_{m=1}^d\|\widetilde{E}_T(\psi_m)\|_{\dot{F}^\sigma}.
\end{equation}
The main result of this subsection is the following proposition.

\newtheorem{Lemmaq7}[Lemmas1]{Proposition}
\begin{Lemmaq7}\label{Lemmaq7} Assume $d\geq 4$. Then, for any $\sigma\in[(d-2)/2,\sigma_0-1]$ and $m=1,\ldots,d$,
\begin{equation}\label{vu20}
\|\mathcal{N}_m(\Psi)\|_{\dot{N}^\sigma[-T,T]}\leq  C\|\widetilde{E}_T(\Psi)\|_{\dot{F}^\sigma}( \|\widetilde{E}_T(\Psi)\|_{\dot{F}^{(d-2)/2}}^2+\|\widetilde{E}_T(\Psi)\|_{\dot{F}^{(d-2)/2}}^4 ).
\end{equation}
\end{Lemmaq7}

The rest of this subsection is concerned with the proof of Proposition \ref{Lemmaq7}. For $\sigma\in[(d-2)/2,\sigma_0-1]$ and $k\in\Z$ let
\begin{equation}\label{go1}
\beta_k(\sigma)=\sum_{m=1}^d\sum_{k'\in\Z}2^{-|k-k'|/10}\cdot (2^{\sigma k'}+2^{(d-2)k'/2})\|\eta_{k'}^{(d)}(\xi)\cdot \mathcal{F}(\widetilde{E}_T(\psi_m))\|_{Z_{k'}}.
\end{equation}
Clearly, $\beta_{k_1}(\sigma)\leq C2^{|k_1-k_2|/10}\beta_{k_2}(\sigma)$ for any $k_1,k_2\in\Z$, and
\begin{equation*}
[\sum_{k\in\Z}\beta_k(\sigma)^2]^{1/2}\leq C\|\widetilde{E}_T(\Psi)\|_{\dot{F}^\sigma}\text{ for any  }\sigma\in[(d-2)/2,\sigma_0-1].
\end{equation*}

For $k\in \Z$ let $P_{k}$ denote the operator defined by the Fourier multiplier $(\xi,\tau)\to\eta_k^{(d)}(\xi)$, and let $P_{\leq k}=\sum_{k'\leq k}P_{k'}$. For $k\in\Z$ and $n\in \Xi_k$ let $\widetilde{P}_{k,n}$ denote the operator defined by the Fourier multiplier $(\xi,\tau)\to\chi_{k,n}(\xi)$. 

\newtheorem{Lemmaq5}[Lemmas1]{Lemma}
\begin{Lemmaq5}\label{Lemmaq5} If $d\geq 4$, $k\in\Z$, $\e'\in\mathbb{S}^{d-1}$, $\sigma\in[(d-2)/2,\sigma_0-1]$, and
\begin{equation}\label{ty1}
F\in\{E_T(A_0),E_T(A_m^2),E_T(\widetilde{\psi}_m\cdot\widetilde{\psi}_l):m,l=1,\ldots,d,\widetilde{\psi}\in\{\psi,\overline{\psi}\}\}
\end{equation}
then
\begin{equation}\label{vu2}
(2^{\sigma k}+2^{(d-2)k/2})\|P_{k}(F)\|_{L^2}\leq C\beta_k(\sigma)\cdot (||\widetilde{E}_T(\Psi)||_{\dot{F}^{(d-2)/2}}+\|\widetilde{E}_T(\Psi)\|_{\dot{F}^{(d-2)/2}}^3),
\end{equation}
and
\begin{equation}\label{vu1}
\|P_{\leq k}(F)\|_{L^{1,\infty}_{\e'}}\leq C2^k(||\widetilde{E}_T(\Psi)||^2_{\dot{F}^{(d-2)/2}}+\|\widetilde{E}_T(\Psi)\|_{\dot{F}^{(d-2)/2}}^4).
\end{equation}
In addition, for $m=1,\ldots,d$,
\begin{equation}\label{vu3}
(2^{\sigma k}+2^{(d-2)k/2})\|P_{k}(E_T(A_m))\|_{L^2}\leq C2^{-k}\beta_k(\sigma)\cdot \|\widetilde{E}_T(\Psi)\|_{\dot{F}^{(d-2)/2}},
\end{equation}
and
\begin{equation}\label{vu10}
\|P_{\leq k}(E_T(A_m))\|_{L^{1,\infty}_{\e'}}\leq C\|\widetilde{E}_T(\Psi)\|_{\dot{F}^{(d-2)/2}}^2.
\end{equation}
\end{Lemmaq5}

The main reason we assume $d\geq 4$ (rather than $d\geq 3$) is to  have a simple proof of  \eqref{vu10}. We defer the proof of Lemma \ref{Lemmaq5} to section \ref{section5}, and complete now the proof of Proposition \ref{Lemmaq7}. For \eqref{vu20} it suffices to prove that
\begin{equation}\label{vu21}
\begin{split}
(2^{\sigma k}+2^{(d-2)k/2})&\|(\tau+|\xi|^2+i)^{-1}\cdot \mathcal{F}(P_k(E_T(\mathcal{N}_m(\Psi))))\|_{Z_k}\\
&\leq C\beta_k(\sigma)\cdot ( \|\widetilde{E}_T(\Psi)\|_{\dot{F}^{(d-2)/2}}^2+\|\widetilde{E}_T(\Psi)\|_{\dot{F}^{(d-2)/2}}^4)
\end{split}
\end{equation}
for any $k\in\Z$. Since $E_T(\mathcal{N}_m(\Psi))$ is  a  sum of terms  of the form $F\cdot \widetilde{E}_T(\psi_m)$ and $E_T(A_l)\cdot\partial_l\widetilde{E}_T(\psi_m)$, where $F$ is as  in \eqref{ty1}, it suffices to prove that 
\begin{equation}\label{ty2}
\begin{split}
&(2^{\sigma k}+2^{(d-2)k/2})\|(\tau+|\xi|^2+i)^{-1}\cdot \mathcal{F}(P_k(F\cdot \widetilde{E}_T(\psi_m)))\|_{Z_k}\\
&+(2^{\sigma k}+2^{(d-2)k/2})\|(\tau+|\xi|^2+i)^{-1}\cdot \mathcal{F}(P_k(E_T(A_l)\cdot \partial_l\widetilde{E}_T(\psi_m)))\|_{Z_k}
\end{split}
\end{equation}
is dominated by the right-hand side of \eqref{vu21} for any  $m,l=1,\ldots,d$. We always estimate the expressions  in \eqref{ty2} using \eqref{im4}. 

We consider first the term $F\cdot \widetilde{E}_T(\psi_m)$, and write $P_k(F\cdot \widetilde{E}_T(\psi_m))$ as
\begin{equation}\label{vu40}
\sum_{|k_1-k|\leq 2}P_k[P_{\leq k-10}(F)\cdot P_{k_1}(\widetilde{E}_T(\psi_m))]+\sum_{k_1\geq k-9}P_k[P_{k_1}(F)\cdot P_{\leq k_1+20}(\widetilde{E}_T(\psi_m))].
\end{equation}
Let $c_\sigma(k)=2^{\sigma k}+2^{(d-2)k/2}$. To  control the term in the first line of  \eqref{ty2} it suffices to prove that for any $v\in I_k^{(d)}$, the quantities
\begin{equation}\label{vu41}
\sum_{|k_1-k|\leq 2}c_\sigma(k)\|\eta_0( |\xi-v|/2^{k-50})(\tau+|\xi|^2+i)^{-1}\mathcal{F}(P_k[P_{\leq k-10}(F)\cdot P_{k_1}(\widetilde{E}_T(\psi_m))])\|_{Z_k}
\end{equation}
and
\begin{equation}\label{vu41.1}
\sum_{k_1\geq  k-9}c_\sigma(k)\|\eta_0( |\xi-v|/2^{k-50})(\tau+|\xi|^2+i)^{-1}\mathcal{F}(P_k[P_{k_1}(F)\cdot P_{\leq k_1+20}(\widetilde{E}_T(\psi_m))])\|_{Z_k}
\end{equation}
are dominated by the right-hand side of \eqref{vu21}. 

To  bound the expression in  \eqref{vu41}, we may assume that $\mathcal{F}(P_{k_1}(\widetilde{E}_T(\psi_m)))$ is supported in $I_{k_1}^{(d)}\times\R\cap
\{(\xi,\tau):|\xi-w| \leq 2^{k_1-50}\}$ for some $w\in
I_{k_1}^{(d)}$.  We use the following simple geometric observation (cf. \cite[Section 8]{IoKe3}): if $\widehat{v},\widehat{w}\in\mathbb{S}^{d-1}$ then there is $\mathbf{e}\in\{\mathbf{e}_1,\ldots,\mathbf{e}_L\}$ such that
\begin{equation}\label{ma2}
\mathbf{e}\cdot\widehat{v}\geq 2^{-5}\text{ and }|\mathbf{e}\cdot\widehat{w}|\geq 2^{-5}.
\end{equation}
We fix $\e$ as in \eqref{ma2} (with $\widehat{v}=v/|v|$ and $\widehat{w}=w/|w|$). Using \eqref{im4}, the expression in \eqref{vu41} is dominated by
\begin{equation*}
\begin{split}
&Cc_\sigma(k)\sum_{|k_1-k|\leq 2}2^{-k/2}\|P_{\leq k-10}(F)\cdot P_{k_1}(\widetilde{E}_T(\psi_m))\|_{L^{1,2}_\e}\\
&\leq Cc_\sigma(k)\sum_{|k_1-k|\leq 2}2^{-k/2}\|P_{\leq k-10}(F)\|_{L^{1,\infty}_\e}\cdot\| P_{k_1}(\widetilde{E}_T(\psi_m))\|_{L^{\infty,2}_\e},
\end{split}
\end{equation*}
which suffices, in view of \eqref{im8} and \eqref{vu1}. 

To bound the expression in \eqref{vu41.1}, we fix $\e\in\{\e_1,\ldots,\e_l\}$ such that $|\e-v/|v||\leq 2^{-100}$ and use \eqref{im4}. The second sum in \eqref{vu41} is dominated by
\begin{equation*}
\begin{split}
&Cc_\sigma(k)\sum_{k_1\geq k-9}2^{-k/2}\|P_k[P_{k_1}(F)\cdot P_{\leq k_1+20}(\widetilde{E}_T(\psi_m))]\|_{L^{1,2}_\e}\\
&\leq Cc_\sigma(k)\sum_{k_1\geq k-9}2^{-k/2}\negmedspace\negmedspace\negmedspace\sum_{n,n'\in\Xi_k\text{ and }|n-n'|\leq C2^k}\|\widetilde{P}_{k,n}P_{k_1}(F)\cdot \widetilde{P}_{k,n'}P_{\leq k_1+20}(\widetilde{E}_T(\psi_m))]\|_{L^{1,2}_\e}\\
&\leq Cc_\sigma(k)\sum_{k_1\geq k-9}2^{-k/2}\|P_{k_1}(F)\|_{L^2}\big[\sum_{n'\in\Xi_k}\|\widetilde{P}_{k,n'}P_{\leq k_1+20}(\widetilde{E}_T(\psi_m))]\|_{L^{2,\infty}_\e}^2\big]^{1/2}\\
&\leq Cc_\sigma(k)\sum_{k_1\geq k-9}2^{-k/2}\frac{\beta_{k_1}(\sigma)\cdot M}{c_\sigma(k_1)}\cdot 2^{k_1/2}2^{-|k_1-k|/4}\|\widetilde{E}_T(\Psi)\|_{\dot{F}^{(d-2)/2}},
\end{split}
\end{equation*}
where $M=(||\widetilde{E}_T(\Psi)||_{\dot{F}^{(d-2)/2}}+\|\widetilde{E}_T(\Psi)\|^3_{\dot{F}^{(d-2)/2}})$, and we used the localized \eqref{im7} and \eqref{vu2} in the last estimate. This suffices since $\beta_{k_1}(\sigma)\leq C2^{|k_1-k|/10}\beta_{k}(\sigma)$ and $d\geq 4$.

We consider now $E_T(A_l)\cdot \partial_l\widetilde{E}_T(\psi_m)$. We write $P_k(E_T(A_l)\cdot\partial_l\widetilde{E}_T(\psi_m))$ as
\begin{equation*}
\begin{split}
&\sum_{|k_1-k|\leq 2}P_k[P_{\leq k-10}(E_T(A_l))\cdot P_{k_1}(\partial_l\widetilde{E}_T(\psi_m))]\\
&+\sum_{k_1\geq k-9}P_k[P_{k_1}(E_T(A_l))\cdot P_{\leq k_1+20}(\partial_l\widetilde{E}_T(\psi_m))],
\end{split}
\end{equation*}
and argue as before, using \eqref{vu10} and \eqref{vu3} instead of \eqref{vu1} and \eqref{vu2}.

\section{Proof of Theorem \ref{Main1}}\label{section4} In this section we assume $d\geq 4$.

\subsection{A priori estimates}\label{section4.1} In this subsection we prove the following:

\newtheorem{Lemmad1}{Proposition}[section]
\begin{Lemmad1}\label{Lemmad1}
Assume that $\sigma_0=d+10$ is as in \eqref{me1}, $T\in[0,1]$ and $s\in C([-T,T]:H^{\infty }_Q)$ is a solution of the initial-value problem
\begin{equation}\label{fi0}
\begin{cases}
&\partial_ts=s\times\Delta s\,\,\text{  on }\,\,\mathbb{R}^d\times[-T,T];\\
&s(0)=s_0.
\end{cases}
\end{equation}
If $\|s_0-Q\|_{\dot{H}^{d/2}}\leq\varepsilon_0\ll 1$ then
\begin{equation}\label{fi20}
\begin{cases}
\sup\limits_{t\in[-T,T]}\|s(t)-Q\|_{\dot{H}^{d/2}}\leq C\|s_0-Q\|_{\dot{H}^{d/2}};\\
\sup\limits_{t\in[-T,T]}\|s(t)\|_{H^{\sigma'}_Q}\leq C( \|s_0\|_{H^{\sigma'}_Q})\text{ for any }\sigma'\in[0,\sigma_0]\cap\Z.
\end{cases}
\end{equation}
\end{Lemmad1}

\begin{proof}[Proof of Proposition \ref{Lemmad1}] We construct $\psi_m,A_m\in C([-T,T]:H^{\infty})$ as in Pro\-position \ref{Lemmag4}. In view of Lemma \ref{Lemmag5}, 
\begin{equation}\label{fi8}
||\psi_{m,0}||_{\dot{H}^{(d-2)/2}}\leq C||s_0-Q||_{\dot{H}^{d/2}}\leq C\varepsilon_0.
\end{equation}
For any $T'\in[0,T]$ we define the functions $E_{T'}(\mathcal{N}_m(\Psi))$ and $\widetilde{E}_{T'}(\psi_m)$ as in \eqref{no6.1} and \eqref{ty4}. Using Propositions \ref{Lemmaqq1} and \ref{Lemmaq3}, for $\sigma\in[(d-2)/2,\sigma_0-1]$ and $T'\in[0,T]$,
\begin{equation}\label{fi1}
\|\widetilde{E}_{T'}(\Psi)\|_{\dot{F}^\sigma}\leq C\cdot ( \sum_{m=1}^d\|\psi_{m,0}\|_{\dot{H}^\sigma\cap\dot{H}^{(d-2)/2}}+\sum_{m=1}^d\|\mathcal{N}_m(\Psi)\|_{\dot{N}^\sigma[-T',T']}).
\end{equation}
In addition, using Lemma \ref{Lemmaq7}, for $\sigma\in[(d-2)/2,\sigma_0-1]$ and $T'\in[0,T]$,
\begin{equation}\label{fi2}
\sum_{m=1}^d\|\mathcal{N}_m(\Psi)\|_{\dot{N}^\sigma[-T',T']}\leq  C\|\widetilde{E}_{T'}(\Psi)\|_{\dot{F}^\sigma}( \|\widetilde{E}_{T'}(\Psi)\|_{\dot{F}^{(d-2)/2}}^2+\|\widetilde{E}_{T'}(\Psi)\|_{\dot{F}^{(d-2)/2}}^4 ).
\end{equation}
The  inequality \eqref{no6.2} shows that the function $L(T')=\sum_{m=1}^d\|\mathcal{N}_m(\Psi)\|_{\dot{N}^\sigma[-T',T']}$ is continuous on the interval $[0,T]$. Also, $L(0)=0$. Thus we can combine \eqref{fi1} and \eqref{fi2} (with $\sigma=(d-2)/2$), together with the smallness of $\|\psi_{m,0}\|_{\dot{H}^{(d-2)/2}}$, to  conclude that
\begin{equation*}
\sum_{m=1}^d\|\mathcal{N}_m(\Psi)\|_{\dot{N}^\sigma[-T',T']}\leq  C\sum_{m=1}^d||\psi_{m,0}||_{\dot{H}^{(d-2)/2}}\text{ for any }T'\in[0,T].
\end{equation*}
Using \eqref{fi1} again, it follows that
\begin{equation}\label{fi3}
\|\widetilde{E}_{T}(\Psi)\|_{\dot{F}^{(d-2)/2}}\leq  C\sum_{m=1}^d||\psi_{m,0}||_{\dot{H}^{(d-2)/2}}\ll 1.
\end{equation}
We combine \eqref{fi1} and \eqref{fi2} again; using \eqref{fi3}, for any $\sigma\in[(d-2)/2,\sigma_0-1]$
\begin{equation}\label{fi4}
\|\widetilde{E}_{T}(\Psi)\|_{\dot{F}^{\sigma}}\leq  C\sum_{m=1}^d||\psi_{m,0}||_{\dot{H}^\sigma\cap\dot{H}^{(d-2)/2}}.
\end{equation}
Using \eqref{im5}, it follows that for any $\sigma\in[(d-2)/2,\sigma_0-1]$
\begin{equation}\label{fi5}
\sum_{m=1}^d\sup_{t\in[-T,T]}\|\psi_m(t)\|_{\dot{H}^\sigma\cap\dot{H}^{(d-2)/2}}\leq C\sum_{m=1}^d||\psi_{m,0}||_{\dot{H}^\sigma\cap\dot{H}^{(d-2)/2}}.
\end{equation}

We use \eqref{fi5} to get a priori estimates on the solution $s$. Using \eqref{fi8} and \eqref{fi5},
\begin{equation}\label{fi6}
\sum_{m=1}^d\sup_{t\in[-T,T]}\|\psi_m(t)\|_{\dot{H}^{(d-2)/2}}\leq C\|s_0-Q\|_{\dot{H}^{d/2}}.
\end{equation}
We define the operators $\nabla^\sigma$, $\sigma\in [-1/2,d/2]$, as in the proof of Lemma \ref{Lemmag5}. Let $p_\sigma=d/(\sigma+1)$. Then, in view of the Sobolev imbedding theorem (recall $d\geq 3$),
\begin{equation}\label{fi7}
\|\nabla^\sigma f\|_{L^{p_\sigma}}\leq C\|\nabla^{\sigma'}f\|_{L^{p_{\sigma'}}}\text{ if }-1/2\leq \sigma\leq\sigma'\leq d/2\text{ and }f\in H^{\sigma_0-1}.
\end{equation}
Let $n_0$ denote the smallest integer $\geq (d-2)/2$. Using \eqref{fi7}, \eqref{t92}, and the definition of the coefficients $A_m$,
\begin{equation}\label{fi9}
\|A_m(t)\|_{\dot{H}^{(d-2)/2}}\leq \|\nabla^{n_0}(A_m(t))\|_{L^{p_{n_0}}}\leq C\|s_0-Q\|_{\dot{H}^{d/2}},
\end{equation}
for any $t\in[-T,T]$ and $m=1,\ldots,d$. 

To prove estimates on the solution $s$, recall the identity \eqref{t49}, 
\begin{equation}\label{fi12}
\begin{cases}
&\partial_ms=\Re(\psi_m)v+\Im(\psi_m)w;\\
&\partial_mv=-\Re(\psi_m)s+A_mw;\\
&\partial_mw=-\Im(\psi_m)s-A_mv.
\end{cases}
\end{equation}
Since $|s|=|v|=|w|\equiv 1$, we use \eqref{fi6}, \eqref{fi9}, and \eqref{fi7} to see that
\begin{equation*}
\sum_{m=1}^d\big[\|\partial_m(s(t))\|_{L^{p_0}}+\|\partial_m(v(t))\|_{L^{p_0}}+\|\partial_m(w(t))\|_{L^{p_0}}\big]\leq C\|s_0-Q\|_{\dot{H}^{d/2}},
\end{equation*}
for any $t\in[-T,T]$. As in the proof of Lemma \ref{Lemmag5}, a simple inductive argument using \eqref{fi12}, \eqref{fi6}, \eqref{fi9}, and \eqref{t92} shows that
\begin{equation}\label{fi14}
\sum_{m=1}^d\big[\|\nabla^{n}\partial_m(s(t))\|_{L^{p_n}}+\|\nabla^{n}\partial_m(v(t))\|_{L^{p_n}}+\|\nabla^{n}\partial_m(w(t))\|_{L^{p_n}}\big]\leq C\|s_0-Q\|_{\dot{H}^{d/2}},
\end{equation}
for any $n\in\Z\cap[0,(d-2)/2]$ and $t\in[-T,T]$. If $d$ is even, this gives
\begin{equation}\label{fi15}
\|s(t)-Q\|_{\dot{H}^{d/2}}\leq C\|s_0-Q\|_{\dot{H}^{d/2}}\text{ for any }t\in[-T,T].
\end{equation}
If $d$ is odd then, using \eqref{fi14} with $n=(d-3)/2$ and \eqref{fi7}, we have
\begin{equation*}
\sum_{m=1}^d\big[\|\nabla^{\sigma}\partial_m(s(t))\|_{L^{p_\sigma}}+\|\nabla^{\sigma}\partial_m(v(t))\|_{L^{p_\sigma}}+\|\nabla^{\sigma}\partial_m(w(t))\|_{L^{p_\sigma}}\big]\leq C\|s_0-Q\|_{\dot{H}^{d/2}}
\end{equation*}
for any $\sigma\in[-1/2,(d-3)/2]$. The bound \eqref{fi15} follows in this case as well, using the Leibniz rule \eqref{Le}.

We show now that for $\sigma'\in[0,\sigma_0]\cap\Z$
\begin{equation}\label{fi30}
\sup_{t\in[-T,T]}\|s(t)\|_{H^{\sigma'}_Q}\leq C( \|s_0\|_{H^{\sigma'}_Q}).
\end{equation}
For this we observe first that we have the conservation law
\begin{equation}\label{conserve}
\|s(t)\|_{H^0_Q}=\|s_0\|_{H^0_Q}\text{ for any }t\in[-T,T],
\end{equation}
which follows by integration by parts from the initial-value problem \eqref{fi0}. Thus, we need to estimate $\|s(t)-Q\|_{\dot{H}^{\sigma'}}$ for $t\in[-T,T]$. Using the first inequality in \eqref{fi20}, we may assume $\sigma'\geq (d+1)/2$. In view of \eqref{t90} and \eqref{fi5}
\begin{equation*}
\sum_{m=1}^d\sup_{t\in[-T,T]}\|\psi_m(t)\|_{\dot{H}^{\sigma'-1}}\leq  C( \|s_0\|_{H^{\sigma'}_Q}).
\end{equation*}
In addition, due to the energy conservation law
\begin{equation*}
\sum_{l=1}^d||\partial_ls(t)||_{L^2}^2=\sum_{l=1}^d||\partial_ls(0)||_{L^2}^2,
\end{equation*}
and the definition $\psi_m=(\partial_ms)\cdot v+i(\partial_ms)\cdot w$, we control $\sup_{t\in[-T,T]}\|\psi_m(t)\|_{L^2}\leq  C( \|s_0\|_{H^{\sigma'}_Q})$. Thus
\begin{equation*}
\sum_{m=1}^d\sup_{t\in[-T,T]}\|\psi_m(t)\|_{H^{\sigma'-1}}\leq  C( \|s_0\|_{H^{\sigma'}_Q}).
\end{equation*}
Using the definition of the coefficients $A_m$, it follows easily that
\begin{equation*}
\sum_{m=1}^d\sup_{t\in[-T,T]}\|A_m(t)\|_{H^{\sigma'-1}}\leq  C( \|s_0\|_{H^{\sigma'}_Q}).
\end{equation*}
We combine the last two inequalities, \eqref{fi12}, and the fact that $|s|=|v|=|w|$; a simple  inductive argument gives $\sup_{t\in[-T,T]}||\partial_ms||_{H^{\sigma'-1}}\leq  C( \|s_0\|_{H^{\sigma'}_Q})$, which completes the proof of \eqref{fi30}.
\end{proof}

\subsection{Existence and uniqueness of solutions}\label{section4.2}

The uniqueness statement in part (a) is proved in \cite[section 2]{IoKe2}: assume $s, s'\in C([T_1,T_2]:H^{\sigma_0}_Q)$ solve the equation $\partial_ts=s\times\Delta_x s$ on $\mathbb{R}^d\times[T_1,T_2]$, and $s(T_1)=s'(T_1)$. Let $q=s'-s$, so
\begin{equation}\label{Sch9}
\begin{cases}
&\partial_tq=(s+q)\times\Delta_x (s+q)-s\times\Delta_x s\,\text{ on }\,\mathbb{R}^d\times[T_1,T_2];\\
&q(T_1)=0.
\end{cases}
\end{equation}
We multiply \eqref{Sch9} by $q(t)$ and integrate by parts over $\mathbb{R}^d$ to obtain
\begin{equation}\label{Sch91}
\begin{split}
\frac{1}{2}\partial_t[\|q(t)\|_{L^2}^2]&=\int_{\mathbb{R}^d}[s(t)\times\Delta_x
q(t)]\cdot q(t)\,dx\\
&\leq C_s(||q(t)||_{L^2}^2+\sum_{l=1}^d||\partial_lq(t)||_{L^2}^2).
\end{split}
\end{equation}
Then we apply $\partial_l$ to \eqref{Sch9}, multiply by $\partial_lq(t)$, add up over $l=1,\ldots,d$, and integrate by parts over  $\mathbb{R}^d$. The result  is
\begin{equation}\label{Sch92}
\begin{split}
\frac{1}{2}\partial_t[\sum_{l=1}^d\|\partial_lq(t)\|_{L^2}^2]&=-\int_{\mathbb{R}^d}[q(t)\times\Delta_x s(t)]\cdot \Delta_xq(t)\,dx\\
&\leq C_s(||q(t)||_{L^2}^2+\sum_{l=1}^d||\partial_lq(t)||_{L^2}^2).
\end{split}
\end{equation}
Using  \eqref{Sch91} and  \eqref{Sch92}, $q\equiv 0$ on $\R^d\times[T_1,T_2]$, as desired.

To construct the global solution, we need the following local existence result:
 
\newtheorem{Lemmad2}[Lemmad1]{Proposition}
\begin{Lemmad2}\label{Lemmad2}
Assume $s_0\in H_Q^{\infty}$. Then there is $T_{\sigma_0}=T( \|s_0\|_{H^{\sigma_0}_Q})>0$ and a solution $s\in C([-T_{\sigma_0},T_{\sigma_0}]:H^{\infty}_Q)$ of the initial-value problem
\begin{equation*}
\begin{cases}
&\partial_ts=s\times\Delta s\,\text{ on }\,\mathbb{R}^d\times[-T_{\sigma_0},T_{\sigma_0}];\\
&s(0)=s_0.
\end{cases}
\end{equation*}
In addition, the time $T_{\sigma_0}$ can be chosen such that
\begin{equation}\label{co1}
\begin{cases}
&\sup\limits_{t\in[-T_{\sigma_0},T_{\sigma_0}]}\|s(t)\|_{H^{\sigma_0}_Q}\leq C(\|s_0\|_{H^{\sigma_0}_Q});\\
&\sup\limits_{t\in[-T_{\sigma_0},T_{\sigma_0}]}\|s(t)\|_{H^\sigma_Q}\leq C(\sigma,\|s_0\|_{H^\sigma_Q})\text{ if  }\sigma\in[\sigma_0,\infty)\cap\Z.
\end{cases}
\end{equation}
\end{Lemmad2}

The local existence Proposition \ref{Lemmad2} is proved, for example, in \cite{KePoStTo}. The bound \eqref{co1} is not stated in this paper, but follows from the key estimate (5.32) in \cite{KePoStTo}. Assuming Proposition \ref{Lemmad2}, by scale invariance, it suffices to construct the solution $s$ in Theorem \ref{Main1} on the time interval $[-1,1]$. In view of Proposition \ref{Lemmad2}, there is $T_{\sigma_0}>0$ and a solution $s$ on the time interval $[-T_{\sigma_0},T_{\sigma_0}]$. Assume the solution $s\in C([T,T]:H^{\infty}_Q)$ is constructed on some time interval $[-T,T]$, $T\leq 1$. In view of Proposition \ref{Lemmad1},
\begin{equation*}
\sup_{t\in[-T,T]}\|s(t)\|_{H^{\sigma_0}_Q}\leq C( \|s_0\|_{H^{\sigma_0}_Q}),
\end{equation*}
uniformly in $T$. Using Proposition \ref{Lemmad2}, the solution $s$ can be extended to the time interval $[-T-T',T+T']$ for some $T'=T'( \|s_0\|_{H^{\sigma_0}})>0$ (which does not depend on $T$). The theorem follows.

\section{Proof of Lemma \ref{Lemmaq5}}\label{section5}

We use the notation in section \ref{section3} and assume in this section that $d\geq 4$. For simplicity of notation, we let $\psi$ denote any of the functions $\widetilde{E}_T(\psi_m)$ or $\overline{\widetilde{E}_T(\psi_m)}$, $m=1,\ldots,d$, $A$ denote any of the functions $A_m$, $m=1,\ldots,d$, and $R$ denote any operator of the form $R_lR_{l'}$, $l,l'=0,1,\ldots,d$, $R_0=I$. With this convention, we show first that for any $k\in\Z$ and $\sigma\in[(d-2)/2,\sigma_0-1]$
\begin{equation}\label{vg1}
(2^{\sigma k}+2^{(d-2)k/2})\|P_k(R(\psi\cdot\psi))\|_{L^2}\leq C\beta_k(\sigma)\|\widetilde{E}_T(\Psi)\|_{\dot{F}^{(d-2)/2}}.
\end{equation}
The left-hand side of \eqref{vg1} is dominated by
\begin{equation*}
\begin{split}
&C(2^{\sigma k}+2^{(d-2)k/2})\sum_{|k_1-k|\leq 2}\sum_{k_2\leq k-4}\|P_{k_1}(\psi)\cdot P_{k_2}(\psi)\|_{L^2}\\
&+C(2^{\sigma k}+2^{(d-2)k/2})\sum_{k_1,k_2\geq k-4,|k_1-k_2|\leq 10}\|P_k(P_{k_1}(\psi)\cdot P_{k_2}(\psi))\|_{L^2}
\end{split}
\end{equation*}
Using \eqref{im8}, we estimate $\|P_{k_1}\psi\|$ in $L^{\infty,2}_\e$ (after suitable localization), and, using the global \eqref{im7}, we estimate $\|P_{k_2}\psi\|$ in $L^{2,\infty}_\e$. The bound \eqref{vg1} follows since $\beta_{k_1}(\sigma)\leq C2^{|k_1-k|/10}\beta_k(\sigma)$. The bounds \eqref{vu2} for $F\in\{E_T(A_0),E_T(\widetilde{\psi}_m\cdot\widetilde{\psi}_l):m,l=1,\ldots,d,\widetilde{\psi}\in\{\psi,\overline{\psi}\}\}$, and \eqref{vu3} clearly follow from \eqref{vg1}. Also, it follows from \eqref{vg1} that
\begin{equation}\label{vg2}
(2^{\sigma k}+2^{(d-2)k/2})\cdot \|P_{k}(A)\|_{L^{\infty,2}_{\e'}}\leq C2^{-k/2}\cdot \beta_k(\sigma)\cdot \|\widetilde{E}_T(\Psi)\|_{\dot{F}^{(d-2)/2}},
\end{equation}
for any $\e'\in\mathbb{S}^{d-1}$.

We prove now that for any $\e'\in\mathbb{S}^{d-1}$
\begin{equation}\label{vg3}
\sum_{k\in\Z}2^{-k}\|P_{k}(R(\psi\cdot\psi))\|_{L^{1,\infty}_{\e'}}\leq C\|\widetilde{E}_T(\Psi)\|^2_{\dot{F}^{(d-2)/2}}.
\end{equation}
For any $k\in\Z$
\begin{equation}\label{vg5}
\begin{split}
\|P_{k}(R(\psi\cdot\psi))\|_{L^{1,\infty}_{\e'}}&\leq C\sum_{|k_1-k|\leq 2}\sum_{k_2\leq k-4}\|P_{k_1}(\psi)\cdot P_{k_2}(\psi)\|_{L^{1,\infty}_{\e'}}\\
&+C\sum_{k_1,k_2\geq k-4,|k_1-k_2|\leq 10}\|P_k(P_{k_1}(\psi)\cdot P_{k_2}(\psi))\|_{L^{1,\infty}_{\e'}}.
\end{split}
\end{equation}
For the first sum in \eqref{vg5}, we use the global \eqref{im7}:
\begin{equation*}
\begin{split}
\sum_{|k_1-k|\leq 2}&\sum_{k_2\leq k-4}\|P_{k_1}(\psi)\cdot P_{k_2}(\psi)\|_{L^{1,\infty}_{\e'}}\\
&\leq C\sum_{|k_1-k|\leq 2}\sum_{k_2\leq k-4}(2^{(d-1)k_1/2}\|P_{k_1}(\psi)\|_{Z_{k_1}})\cdot (2^{(d-1)k_2/2}\|P_{k_2}(\psi)\|_{Z_{k_2}})\\
&\leq C2^k\beta_k((d-2)/2)^2.
\end{split}
\end{equation*}
For the second sum, we use the localized \eqref{im7} and the assumption $d\geq 4$:
\begin{equation*}
\begin{split}
\|P_k(P_{k_1}\psi\cdot &P_{k_2}\psi)\|_{L^{1,\infty}_{\e'}}\leq\,C\sum_{n,n'\in\Xi_k\text{ and }|n-n'|\leq C2^k}\|\widetilde{P}_{k,n}P_{k_1}(\psi)\cdot \widetilde{P}_{k,n'}P_{k_2}(\psi)\|_{L^{1,\infty}_{\e'}}\\
&\leq C\big[\sum_{n\in\Xi_k}\|\widetilde{P}_{k,n}P_{k_1}(\psi)\|^2_{L^{2,\infty}_{\e'}}\big]^{1/2}\cdot\big[\sum_{n'\in\Xi_k}\|\widetilde{P}_{k,n}P_{k_2}(\psi)\|^2_{L^{2,\infty}_{\e'}}\big]^{1/2}\\
&\leq C2^{-3|k_1-k|/2}\cdot (2^{(d-1)k_1/2}\|P_{k_1}(\psi)\|_{Z_{k_1}})\cdot (2^{(d-1)k_2/2}\|P_{k_2}(\psi)\|_{Z_{k_2}})\\
&\leq C2^k2^{-|k_1-k|/4}\beta_k((d-2)/2)^2. 
\end{split}
\end{equation*}
The bound \eqref{vg3} follows from \eqref{vg5} and the last two estimates. The bounds \eqref{vu1} for $F\in\{E_T(A_0),E_T(\widetilde{\psi}_m\cdot\widetilde{\psi}_l):m,l=1,\ldots,d,\widetilde{\psi}\in\{\psi,\overline{\psi}\}\}$, and \eqref{vu10} clearly follow from \eqref{vg3}.

It remains to prove the bounds \eqref{vu2} and \eqref{vu1} for $F=E_T(A_m^2)$. We will need the following technical lemma:

\newtheorem{Lemmaz5}{Lemma}[section]
\begin{Lemmaz5}\label{Lemmaz5}
If $k\in\Z$, $k'\in(-\infty, k+10d]\cap\Z$, and $\e'\in\mathbb{S}^{d-1}$ then
\begin{equation}\label{jj1}
\big[\sum_{n\in\Xi_{k'}}\|\widetilde{P}_{k',n}P_{k}(A)\|^2_{L^{2,\infty}_{\e'}}\big]^{1/2}\leq C2^{k/2}2^{-3|k-k'|/4}\|\widetilde{E}_T(\Psi)\|_{\dot{F}^{(d-2)/2}}^2.
\end{equation}
\end{Lemmaz5}

Assuming Lemma \ref{Lemmaz5}, for \eqref{vu2} it suffices to prove that
\begin{equation}\label{vg9}
(2^{\sigma k}+2^{(d-2)k/2})\|P_k(A\cdot A)\|_{L^2}\leq C\beta_k(\sigma)\|\widetilde{E}_T(\Psi)\|^3_{\dot{F}^{(d-2)/2}}.
\end{equation}
The proof of \eqref{vg9} is similar to the proof of \eqref{vg1}, using the $L^{\infty,2}_{\e'}$ estimate in \eqref{vg2} and the global (that is $k'=k$) $L^{2,\infty}_{\e'}$ estimate in \eqref{jj1}. For \eqref{vu1} it suffices to prove that
\begin{equation}\label{vg10}
\|P_{k}(A\cdot A)\|_{L^{1,\infty}_{\e'}}\leq C2^k\|\widetilde{E}_T(\Psi)\|^4_{\dot{F}^{(d-2)/2}},
\end{equation}
for any $k\in\Z$ and $\e'\in\mathbb{S}^{d-1}$. The  proof of \eqref{vg10} is similar to the proof of \eqref{vg3}, using the localized $L^{2,\infty}_{\e'}$ estimate in \eqref{jj1}.

\begin{proof}[Proof of Lemma \ref{Lemmaz5}] In view of the definitions, we may assume $k'\leq k-10d$ and it suffices to prove that
\begin{equation}\label{jj2}
\big[\sum_{n\in\Xi_{k'}}\|\widetilde{P}_{k',n}P_{k}(\psi\cdot \psi)\|^2_{L^{2,\infty}_{\e'}}\big]^{1/2}\leq C2^{3k/2}2^{-3|k-k'|/4}\|\widetilde{E}_T(\Psi)\|^2_{\dot{F}^{(d-2)/2}}.
\end{equation}
We will use the following bound: if $k\in\Z$, $k'\in(-\infty,k+10d]\cap\Z$, and $f\in Z_k$ then
\begin{equation}\label{im6}
\big[\|\sum_{n\in\Xi_{k'}}\mathcal{F}^{-1}(\chi_{k',n}(\xi)\cdot \widetilde{f})\|_{L^\infty_{x,t}}^2\big]^{1/2}
\leq C2^{dk/2}\cdot 2^{-d|k-k'|/2}(1+|k-k'| )\cdot\|f\|_{Z_k},
\end{equation}
where $\mathcal{F}^{-1}(\widetilde{f})\in\{\mathcal{F}^{-1}(f),\overline{\mathcal{F}^{-1}(f)}\}$. For $k-k'\leq C$ this follows directly from \eqref{im5} and the Sobolev imbedding theorem. For $k-k'\geq C$, the bound \eqref{im6} follows by analyzing the cases $f\in X_k$ and $f\in Y_k^\e$ (see Lemma 4.1 in \cite{IoKe3} for a similar proof).

The left-hand side of \eqref{jj2} is dominated by
\begin{equation}\label{jj3}
\begin{split}
&C\sum_{|k_1-k|\leq 2}\sum_{k_2\leq k'}\big[\sum_{n\in\Xi_{k'}}\|\widetilde{P}_{k',n}P_{k}(P_{k_1}(\psi)\cdot P_{k_2}(\psi))\|^2_{L^{2,\infty}_{\e'}}\big]^{1/2}\\
&+C\sum_{|k_1-k|\leq 2}\sum_{k'\leq k_2\leq k-4}\big[\sum_{n\in\Xi_{k'}}\|\widetilde{P}_{k',n}P_{k}(P_{k_1}(\psi)\cdot P_{k_2}(\psi))\|^2_{L^{2,\infty}_{\e'}}\big]^{1/2}\\
&+C\sum_{k_1,k_2\geq k-4,\,|k_1-k_2|\leq 10}\big[\sum_{n\in\Xi_{k'}}\|\widetilde{P}_{k',n}P_{k}(P_{k_1}(\psi)\cdot P_{k_2}(\psi))\|^2_{L^{2,\infty}_{\e'}}\big]^{1/2}.
\end{split}
\end{equation}
We use the $L^\infty_{x,t}$ estimate \eqref{im6} on the lower frequency term and the localized $L^{2,\infty}_{\e'}$ estimate \eqref{im7} on the higher frequency term. The first sum in \eqref{jj3}  is dominated by
\begin{equation*}
C\sum_{|k_1-k|\leq 2}\sum_{k_2\leq k'}(2^{k_1/2}2^{-3|k-k'|/4}\|\widetilde{E}_T(\Psi)\|_{\dot{F}^{(d-2)/2}})\cdot (2^{k_2}\|\widetilde{E}_T(\Psi)\|^2_{\dot{F}^{(d-2)/2}}),
\end{equation*}
which suffices for \eqref{jj2}. The second sum in \eqref{jj3}  is dominated by
\begin{equation*}
\begin{split}
&C\sum_{|k_1-k|\leq 2}\sum_{k'\leq k_2\leq k-4}\big[\sum_{n\in\Xi_{k'}}\|\widetilde{P}_{k',n}P_{k_1}(\psi)\|^2_{L^{2,\infty}_{\e'}}\big]^{1/2}\cdot \big[\sum_{n\in\Xi_{k'}}\|\widetilde{P}_{k',n}P_{k_2}(\psi)\|_{L^\infty}\big]\\
&\leq C\sum_{|k_1-k|\leq 2}\sum_{k'\leq k_2\leq k-4}(2^{k/2}2^{-7|k-k'|/8}\|\widetilde{E}_T(\Psi)\|_{\dot{F}^{(d-2)/2}})\negmedspace\cdot\negmedspace (2^{k_2}|k-k'|\|\widetilde{E}_T(\Psi)\|_{\dot{F}^{(d-2)/2}})
\end{split}
\end{equation*}
which suffices for \eqref{jj2}. The third sum in \eqref{jj3}  is dominated by
\begin{equation*}
\begin{split}
&C2^{d|k-k'|/2}\negmedspace\negmedspace\negmedspace\sum_{k_1,k_2\geq k-4,\,|k_1-k_2|\leq 10}\negmedspace\big[\sum_{n\in\Xi_{k'}}\|\widetilde{P}_{k',n}P_{k_1}(\psi)\|^2_{L^{2,\infty}_{\e'}}\big]^{1/2}\negmedspace\cdot \negmedspace\big[\sum_{n\in\Xi_{k'}}\|\widetilde{P}_{k',n}P_{k_2}(\psi)\|_{L^\infty}^2\big]^{1/2}\\
&\leq C2^{d|k-k'|/2}\negmedspace\negmedspace\sum_{k_1,k_2\geq k-4,\,|k_1-k_2|\leq 10}\negmedspace2^{3k_1/2}\|\widetilde{E}_T(\Psi)\|_{\dot{F}^{(d-2)/2}}^2\cdot 2^{-(d-1)|k_1-k'|}(1+|k_1-k'|)^2,
\end{split}
\end{equation*}
which suffices for \eqref{jj2} since $d\geq 4$. This completes the proof of Lemma \ref{Lemmaz5}.
\end{proof}

\end{document}